\documentclass[11pt]{amsart}

\usepackage{amsmath,amssymb,amsthm}
\usepackage[T1]{fontenc}
\usepackage{lmodern}
\usepackage{microtype}
\usepackage[dvipsnames]{xcolor}
\usepackage{hyperref}
\definecolor{darkblue}{rgb}{0,0,0.7}
\hypersetup{colorlinks=true,citecolor=darkblue,linkcolor=darkblue,
 urlcolor=darkblue,
 pdftitle={Kronecker coefficients of large degree and Littlewood--Richardson inequalities},
 pdfauthor={Zachary Slonim}}
\newcommand{\defn}[1]{\emph{\textcolor{darkblue}{#1}}}

\makeatletter
\def\th@plain{\thm@notefont{}\itshape}
\def\th@definition{\thm@notefont{}\normalfont}
\makeatother
\newtheorem{thm}{Theorem}[section]
\newtheorem{lemma}[thm]{Lemma}
\newtheorem{prop}[thm]{Proposition}
\newtheorem{cor}[thm]{Corollary}
\theoremstyle{definition}
\newtheorem{remark}[thm]{Remark}
\newtheorem{example}[thm]{Example}
\numberwithin{equation}{section}

\newcommand{\Pl}{\operatorname{Pl}}
\newcommand{\supp}{\operatorname{supp}}
\newcommand{\Ind}{\operatorname{Ind}}
\newcommand{\Res}{\operatorname{Res}}
\newcommand{\sgn}{\operatorname{sgn}}
\newcommand{\tr}{\operatorname{tr}}
\newcommand{\T}{\mathsf T}
\newcommand{\HH}{\mathcal H}
\newcommand{\CC}{\mathbb C}
\newcommand{\EE}{\mathbb E}
\newcommand{\one}{\mathbf 1}
\newcommand{\dGamma}{\mathrm d\Gamma}
\newcommand{\ts}{\hskip0.03cm}

\title[Kronecker coefficients of large degree]
{Kronecker coefficients of large degree\\
and Littlewood--Richardson inequalities}
\author[Zachary Slonim]{Zachary Slonim$^\star$}
\date{}
\thanks{November 5, 2026.}
\thanks{\thinspace ${\hspace{-.45ex}}^\star$Department of Mathematics,
UCLA, Los Angeles, CA 90095, USA. Email: \texttt{zachslonim@math.ucla.edu}}

\begin{document}
\begin{abstract}
We prove that Kronecker coefficients are uniformly asymptotic to
$f^\lambda f^\mu f^\nu/n!$ when all three irreducible dimensions are at
least $\sqrt{n!}\,e^{-A\sqrt n}$, for any fixed $A>0$. This strengthens a
conjecture of Pak, Panova and Yeliussizov. For triple staircase and triple
square partitions, we obtain relative asymptotics with respective
limiting factors $e^{-1/24}$ and $e^{-1/3}$. For Littlewood--Richardson (LR)
coefficients, we bound the sum of the three dimension deficits by the
LR coefficient deficit times $O(\log n)$, proving their Conjecture~5.5 and
extending it to arbitrary fixed size proportions. We also construct
vanishing LR coefficients with three Plancherel dimension sequences and
determine the exact stabilization threshold of the largest LR coefficient
with one input size fixed.
\end{abstract}
\maketitle

\section{Introduction}\label{sec:intro}

\subsection{Foreword}
Kronecker and Littlewood--Richardson coefficients describe tensor
product and induction for symmetric groups, or equivalently the
internal and ordinary products of Schur functions
\cite{LR34,Mur38,Mac95}. They connect representation theory with
tableaux, geometry and complexity \cite{Ful00,BMS15,Pan24}.
Two questions guide their study: which multiplicities are positive,
and how large can they be?

Irreducible dimensions provide a benchmark for size.
For $\lambda\vdash n$, $f^\lambda$ counts standard Young tableaux,
and $(f^\lambda)^2/n!$ is the Robinson--Schensted shape distribution
of a uniform permutation \cite{Sag01}. Dimension bounds and the
VKLS limit shape describe large-degree partitions
\cite{LS77,VK81,VK85}; the precise leading constant for the maximal
dimension is now known \cite{AE26}. Pak, Panova and Yeliussizov
relate nearly maximal multiplicities to large dimensions
\cite{PPY19}, with analogous questions for finite groups in
\cite{PPY19b}. Our main theorem gives the converse: three large dimensions
determine a Kronecker coefficient up to a uniform vanishing
relative error.

Explicit shapes suggest finer asymptotic questions.
Rectangular coefficients connect Kronecker products with Gaussian
binomials and partition counts \cite{PP14,PP17,MePP20}; contingency
tables and Durfee size give other growth bounds \cite{PP20,PP23}.
For staircases, Saxl's conjecture asks whether the tensor square
contains every irreducible. The study of its constituents
\cite{PPV16,Ike15,LuS17}, and a recent proof preprint \cite{OAI26},
motivate the quantitative question of how many copies occur.
We determine the leading term for triple staircase and triple
square coefficients, sharpening the growth conjecture in
\cite{PP23}. Character methods underlie these estimates \cite{Bia98,FS11,MS22}.

For LR coefficients, tableau and polyhedral models connect
positivity with saturation and Horn's eigenvalue problem
\cite{KT99,PV05}. Their sizes also enter skew-tableau enumeration,
weighted lozenge tilings and large-deviation theory
\cite{MPP18,MPT22,BGH22}. Extremal questions bring these viewpoints
together: Schur-positivity inequalities constrain maximizing
diagrams \cite{LPP07,PS26}, while dimension estimates constrain
nearly maximal coefficients \cite{PPY19}. We quantify the latter
relation through a dimension-deficit inequality, and show that all
three dimensions can be Plancherel large even when an LR
coefficient vanishes. We also determine when the largest LR coefficient first stabilizes
with one input size fixed.

\newpage
\subsection{Kronecker coefficients of large degree}
For $\lambda\vdash n$, let $V^\lambda$ be the irreducible complex
$S_n$-module indexed by $\lambda$, and put $f^\lambda:=\dim V^\lambda$.
The \defn{Kronecker coefficient} $g(\lambda,\mu,\nu)$ is the multiplicity
of $V^\nu$ in $V^\lambda\otimes V^\mu$. Define
\begin{equation}\label{eq:intro-RH}
 R_n(\lambda,\mu,\nu)
   :=\frac{n!\,g(\lambda,\mu,\nu)}{f^\lambda f^\mu f^\nu},
 \qquad
 \HH_A(n):=\{\lambda\vdash n:
       f^\lambda\ge\sqrt{n!}\,e^{-A\sqrt n}\}.
\end{equation}

\begin{thm}[Uniform asymptotic]\label{thm:uniform}
For every fixed $A>0$, as $n\to\infty$,
\begin{equation}\label{eq:uniform}
 \max_{\lambda,\mu,\nu\ts\in\ts \HH_A(n)}
       |R_n(\lambda,\mu,\nu)-1|\longrightarrow0.
\end{equation}
In particular, all the coefficients in this maximum are positive for
sufficiently large $n$. An empty maximum of nonnegative errors is
interpreted as zero.
\end{thm}

This strengthens \cite[Conjecture~3.10]{PPY19}: it identifies the leading
term uniformly, whereas the conjecture asks for a lower bound on the
exponential scale. We prove Theorem~\ref{thm:uniform} in
Section~\ref{sec:kronecker}. Its dimension-scale consequence is the
following.

Following \cite[Definition~1.2]{PPY19}, a sequence
$\lambda^{(n)}\vdash n$ is a \defn{Plancherel sequence} if
$f^{\lambda^{(n)}}\ge\sqrt{n!}\,e^{-A\sqrt n}$ for some fixed $A>0$.
This is a deterministic dimension condition.

\begin{cor}[Three Plancherel sequences]\label{cor:ppy310}
Let $\lambda^{(n)},\mu^{(n)},\nu^{(n)}\vdash n$ be three Plancherel
sequences. Then
\begin{equation}\label{eq:ppy310}
 g(\lambda^{(n)},\mu^{(n)},\nu^{(n)})
  =(1+o(1))\frac{f^{\lambda^{(n)}}f^{\mu^{(n)}}f^{\nu^{(n)}}}{n!}
  =\sqrt{n!}\,e^{-O(\sqrt n)}.
\end{equation}
\end{cor}

Thus \cite[Conjecture~3.10]{PPY19} holds with its loss $O(n)$ improved
to $O(\sqrt n)$. More precisely, if
$K(n):=\max_{\alpha,\beta,\gamma\ts\vdash\ts n}g(\alpha,\beta,\gamma)$, then
every triple of Plancherel sequences satisfies
\[
 e^{-O(\sqrt n)}\le
 \frac{g(\lambda^{(n)},\mu^{(n)},\nu^{(n)})}{K(n)}\le1.
\]
Thus every such coefficient is within a multiplicative factor
$e^{O(\sqrt n)}$ of the largest Kronecker coefficient.

\subsection{Staircase and square partitions}
Let
\[
 \rho_k:=(k-1,k-2,\ldots,1),\qquad
 \delta_k:=(k^k),\qquad
 n_\rho:=\binom k2,\qquad n_\delta:=k^2.
\]
For these explicit families, the relative asymptotic has a constant
factor arising from permutations whose nontrivial cycles all have
length three.

\begin{thm}[Triple staircases and triple squares]\label{thm:shapes}
As $k\to\infty$,
\begin{align}
 g(\rho_k,\rho_k,\rho_k)
   &=\bigl(e^{-1/24}+o(1)\bigr)
           \frac{(f^{\rho_k})^3}{n_\rho!},
       \label{eq:staircase-asymptotic}\\
 g(\delta_k,\delta_k,\delta_k)
   &=\bigl(e^{-1/3}+o(1)\bigr)
           \frac{(f^{\delta_k})^3}{n_\delta!}.
       \label{eq:square-asymptotic}
\end{align}
\end{thm}

The proof is in Section~\ref{sec:shapes}. These families have
dimension deficit proportional to their size, as computed in
Proposition~\ref{prop:shape-dimensions}. Thus the theorem treats a
different dimension range from Theorem~\ref{thm:uniform}. For either
family, writing $n$ for its size, the resulting coefficient is
$\sqrt{n!}\,e^{-\Theta(n)}$. In particular, this proves
\cite[Conjecture~6.2]{PP23}, with the stronger relative asymptotics
\eqref{eq:staircase-asymptotic} and \eqref{eq:square-asymptotic}.

\subsection{Littlewood--Richardson coefficients}
For $\mu\vdash k$, $\nu\vdash l$ and $\lambda\vdash n$, with
$n:=k+l$, the \defn{Littlewood--Richardson coefficient}
$c^\lambda_{\mu\nu}$ is the multiplicity of $V^\lambda$ in
$\Ind_{S_k\times S_l}^{S_n}(V^\mu\boxtimes V^\nu)$.
For a partition $\rho\vdash r$, put
\[
 \Delta(\rho):=\log\frac{\sqrt{r!}}{f^\rho}.
\]
If $c^\lambda_{\mu\nu}>0$, put
\begin{equation}\label{eq:intro-delta}
 \delta(\lambda;\mu,\nu)
       :=\log\frac{\sqrt{\binom nk}}{c^\lambda_{\mu\nu}}.
\end{equation}
Both quantities are nonnegative; the universal coefficient bound is
recalled in Proposition~\ref{prop:mackey} below. All logarithms are
natural unless a base is specified.

\begin{thm}[Dimension deficits]\label{thm:deficits}
For every fixed $0<\varepsilon\le1/2$, there are constants
$a_\varepsilon,b_\varepsilon>0$ such that, whenever
$n=k+l\ge2$, $k,l\ge\varepsilon n$, and $c^\lambda_{\mu\nu}>0$,
\begin{equation}\label{eq:deficits}
 \Delta(\lambda)+\Delta(\mu)+\Delta(\nu)
   \le a_\varepsilon+
       b_\varepsilon\,\delta(\lambda;\mu,\nu)\log n.
\end{equation}
\end{thm}

This proves \cite[Conjecture~5.5]{PPY19}, strengthening it to all three
dimensions and to size proportions bounded away from zero and one.
The proof is in Section~\ref{sec:deficits}.
Indeed, when $k=l=n/2$, exponentiation gives universal $a,u>0$ with
\begin{equation}\label{eq:ppy55}
 \frac{f^\lambda}{\sqrt{n!}}
    \ge a\left(\frac{c^\lambda_{\mu\nu}}
                         {\sqrt{\binom n{n/2}}}\right)^{u\log n}.
\end{equation}
This is the conjectured inequality, also valid when the coefficient is
zero. Theorem~\ref{thm:deficits} also gives the extension of
\cite[Theorem~1.7]{PPY19} to arbitrary fixed proportions requested in
\cite[\S7.7]{PPY19}; see Corollary~\ref{cor:proportions}.

\begin{thm}[Vanishing with three large dimensions]\label{thm:vanishing}
For every fixed $A>0$, there is $A'>0$ such that, for every sufficiently
large even $n$ and every $\mu,\nu\in\HH_A(n/2)$, there is
$\lambda\in\HH_{A'}(n)$ satisfying
\begin{equation}\label{eq:vanishing}
 c^\lambda_{\mu\nu}=0.
\end{equation}
\end{thm}

This proves \cite[Conjecture~4.16]{PPY19} with the extended-real
convention $\log0:=-\infty$, and strengthens
\cite[Proposition~4.17]{PPY19} by reducing its output dimension deficit
from $O(\sqrt n\log n)$ to $O(\sqrt n)$. We prove
Theorem~\ref{thm:vanishing} in Section~\ref{sec:vanishing}.
More explicitly, one may take
three Plancherel sequences, indexed by even $n$, for which
\[
 \frac1{\sqrt n}\left(\frac n2-
            \log_2 c^{\lambda^{(n)}}_{\mu^{(n)}\nu^{(n)}}\right)
        \longrightarrow+\infty.
\]
The role of the convention at zero is discussed in
Remark~\ref{rem:zero}.

For integers $n\ge k\ge0$, define
\begin{equation}\label{eq:intro-maxima}
 D(k):=\max_{\mu\ts\vdash\ts k}f^\mu,
 \qquad
 C(n,k):=\max_{\substack{\lambda\ts\vdash\ts n,\ \mu\ts\vdash\ts k\\
                         \nu\ts\vdash\ts n-k}}c^\lambda_{\mu\nu},
 \qquad N_k:=\binom{k+1}{2}.
\end{equation}
We allow the empty partition and set $D(0):=1$.

\begin{thm}[Exact stabilization threshold]\label{thm:stabilization}
For $k\ge3$ and $n\ge k$,
\begin{equation}\label{eq:stabilization}
 C(n,k)=D(k)\quad\Longleftrightarrow\quad n\ge N_k.
\end{equation}
For $k=0,1,2$, one has $C(n,k)=D(k)=1$ for every $n\ge k$.
\end{thm}

This proves \cite[Remark~4.24]{PPY19} for $k\ge3$, adding necessity
to \cite[Theorem~4.23]{PPY19}. For $k=2$, the first occurrence is
$n=2<N_2=3$. The proof is in Section~\ref{sec:stabilization}.

\newpage
\subsection{Proof idea}\label{subsec:proof-idea}
For Theorem~\ref{thm:uniform}, the starting point is that the identity
permutation contributes exactly $1$ to $R_n$. To show that
all other permutations contribute $o(1)$ uniformly, we prove that the
cubic absolute sum of a normalized irreducible character tends to zero.
A large dimension first bounds the first row and column by $O(\sqrt n)$.
Specializations of the Cauchy identities then force cancellation in the
content formulas for transpositions and $3$-cycles. We retain the
transposition cancellation in the Stanley--F\'eray character formula,
using the factorization estimates of F\'eray and \'{S}niady
\cite[Theorem~2, Proposition~12 and Lemma~13]{FS11}.
This controls permutations of small support. Orthogonality
and conjugacy-class sizes handle the remaining permutations. H\"older's
inequality gives Theorem~\ref{thm:uniform}.

For Theorem~\ref{thm:shapes}, the character on one $3$-cycle determines
the contributions of any fixed number of disjoint $3$-cycles through
Lemma~\ref{lem:three-cycles}. Summing these contributions gives the
two exponential constants. Uniform bounds for all supports justify
the summation: for staircases, we combine Murnaghan--Nakayama with a
dimension-ratio estimate; for squares, an exterior-power formula gives
an extra decay factor for every even cycle.

For Theorem~\ref{thm:deficits}, we regard the squared character as a
central probability distribution on a symmetric group. A large dimension
deficit forces most of this distribution away from permutations of small
support. Such permutations rarely preserve a subset of proportional
size, which bounds a restriction norm. Mackey's formula expresses an
induction norm in terms of two such restrictions; concentration of the
overlap gives the required inequality.

For Theorem~\ref{thm:vanishing}, we force LR vanishing by adding a
row of length $O(\sqrt n)$ to a large-dimensional diagram and use the
hook-length formula to control the cost.

For Theorem~\ref{thm:stabilization}, we show that
equality in $c^\lambda_{\mu\nu}\le f^\mu$, when $f^\mu>1$, forces
$\lambda/\nu$ to have at most one cell in each row and column. The
smallest outer partition containing such a skew diagram with $k$ cells
has size $N_k$.

\subsection{Structure of the paper}
Section~\ref{sec:notation} fixes the notation. In
Section~\ref{sec:background}, we record the dimension estimates,
character identities and restriction constructions used in the proofs.
Section~\ref{sec:kronecker} proves Theorem~\ref{thm:uniform} and
Corollary~\ref{cor:ppy310}. Section~\ref{sec:shapes} proves the
staircase and square asymptotics of Theorem~\ref{thm:shapes}.
Sections~\ref{sec:deficits},
\ref{sec:vanishing} and~\ref{sec:stabilization} prove
Theorems~\ref{thm:deficits}, \ref{thm:vanishing} and
\ref{thm:stabilization}, respectively.

\section*{Declaration of AI use}
All mathematical statements and proofs were found by ChatGPT 6 Pro.
We do not claim any credit for the results. This manuscript was written
almost entirely by ChatGPT-6.1 Sol Max in an iterative process directed
by the author in the hopes of making the proofs more understandable.
No independent verification of the technical proofs has been done.

\section{Basic definitions and notation}\label{sec:notation}

\subsection{Partitions and tableaux}
We write $\lambda\vdash n$ for a partition of $n$, $|\lambda|$ for its
size, $\ell(\lambda)$ for its number of nonzero parts, and
$\lambda'$ for its conjugate, whose $j$th part is
$\lambda'_j:=|\{i:\lambda_i\ge j\}|$. The \defn{Young diagram} of
$\lambda$ is the left-justified array with $\lambda_i$ cells in row $i$.
Rows are numbered from top to bottom and columns from left to right;
conjugation interchanges the rows and columns. For a cell $b=(i,j)$,
its \defn{content} and \defn{hook length}
are
\[
 c(b):=j-i,
 \qquad h_\lambda(i,j):=\lambda_i-j+\lambda'_j-i+1.
\]
For $\nu\subseteq\lambda$, the skew diagram $\lambda/\nu$ is a
\defn{rook strip} if it has at most one cell in each row and each column.
A \defn{standard Young tableau} of shape $\lambda/\nu$, with
$r:=|\lambda/\nu|$, is a filling of its cells with $1,\ldots,r$, each used
exactly once, whose entries increase from left to right along each row
and from top to bottom along each column. For a straight shape
$\lambda\vdash n$, the number of standard Young tableaux of shape
$\lambda$ equals the irreducible dimension $f^\lambda:=\dim V^\lambda$.

We write $s_\lambda$ and $s_{\lambda/\nu}$ for Schur and skew Schur
functions, respectively. The LR coefficients defined by induction are
also the coefficients in the following Schur expansions:
\begin{equation}\label{eq:Schur-LR}
 s_\mu s_\nu=\sum_\lambda c^\lambda_{\mu\nu}s_\lambda,
 \qquad
 s_{\lambda/\nu}=\sum_\mu c^\lambda_{\mu\nu}s_\mu.
\end{equation}
For the Schur expansions, see \cite[Theorems~4.7.4, 4.9.2 and~4.9.4]{Sag01};
for the tableau and representation-theoretic constructions, see
\cite[\S\S2.1--2.3]{Sag01} and \cite[Chapter~I, \S7]{Mac95}. All representations below are
over $\CC$, so restriction and induction are semisimple.

\subsection{Characters and dimension deficits}\label{subsec:characters}
Let $\chi^\lambda$ be the character of $V^\lambda$ and define
\begin{equation}\label{eq:normalized-character}
 r_\lambda(\pi):=\frac{\chi^\lambda(\pi)}{f^\lambda},
 \qquad
 \Delta(\lambda):=\log\frac{\sqrt{|\lambda|!}}{f^\lambda}.
\end{equation}
Thus $r_\lambda(e)=1$ and $|r_\lambda(\pi)|\le1$. We put
$f^\varnothing:=1$ and $\Delta(\varnothing):=0$. The Plancherel measure is
\[
 \Pl_n(\lambda):=\frac{(f^\lambda)^2}{n!}=e^{-2\Delta(\lambda)}.
\]
In particular, $\Delta(\lambda)\ge0$ and
$\HH_A(n)=\{\lambda\vdash n:\Delta(\lambda)\le A\sqrt n\}$.
A partition of $n$ is \defn{$B$-balanced} if both $\lambda_1$ and
$\ell(\lambda)$ are at most $B\sqrt n$.

For $\pi\in S_n$, write $s(\pi):=|\supp\pi|$, and let $|\pi|$ denote
its minimum length as a product of transpositions. If its nontrivial
cycles number $q$, then $|\pi|=s(\pi)-q$. Write
\[
 (n)_s:=n(n-1)\cdots(n-s+1),
 \qquad \Sigma_\lambda(\pi):=(n)_{s(\pi)}r_\lambda(\pi).
\]
In the argument of $\Sigma_\lambda$, fixed points are omitted; for
example, $\Sigma_\lambda(3^q)$ denotes $q$ disjoint $3$-cycles.
We set $\Sigma_\lambda(\varnothing):=1$.

All asymptotics refer to $n\to\infty$. A subscript on $O(\cdot)$
specifies the permitted dependence of its constant. Every estimate over
$\HH_A(n)$ is uniform over that set. Maxima of nonnegative errors over
an empty set are assigned the value zero.

\section{Background results and constructions}\label{sec:background}

We collect the ingredients used in the proofs. Dimension estimates
and row removal turn large degree into balancedness, while
orthogonality and class sizes control characters on large supports.
For small supports, we introduce content specializations and the
Stanley--F\'eray coloring formula. Finally, restriction and induction
norms connect LR coefficients with the dimensions of their inputs
and output.

\subsection{Dimensions and balanced diagrams}
The hook-length formula \cite[Theorem~3.10.2]{Sag01} and the
hook-content formula \cite[Chapter~I, \S3, Example~4]{Mac95} give
\begin{equation}\label{eq:hooks}
 f^\lambda=\frac{n!}{\prod_{b\ts\in\ts \lambda}h_\lambda(b)},
 \qquad
 s_\lambda(1^z)=\frac{f^\lambda}{n!}
                       \prod_{b\ts\in\ts \lambda}(z+c(b)).
\end{equation}
Here $s_\lambda(1^z)$ denotes the polynomial specialization
$p_j\mapsto z$ for every power sum $p_j$; for positive integral $z$,
it is evaluation at $z$ variables equal to one.

Recall $D(n):=\max_{\lambda\ts\vdash\ts n}f^\lambda$, and let
$p(n)$ denote the number of partitions of $n$. We shall use
\begin{equation}\label{eq:partition-bound}
 p(n)\le e^{\pi\sqrt{2n/3}},
 \qquad
 D(n)\ge\sqrt{n!}\,e^{-a_0\sqrt n},
 \qquad a_0:=\frac\pi{\sqrt6}.
\end{equation}
For completeness, Euler's generating function and $e^x-1\ge x$ give,
for $t>0$,
\[
 p(n)e^{-tn}\le\prod_{j\ts\ge\ts 1}(1-e^{-tj})^{-1},
 \qquad
 \log\prod_{j\ts\ge\ts 1}(1-e^{-tj})^{-1}
   =\sum_{r\ts\ge\ts 1}\frac1{r(e^{tr}-1)}\le\frac{\pi^2}{6t}.
\]
Take $t:=\pi/\sqrt{6n}$. The second assertion of
\eqref{eq:partition-bound} follows from $n!\le p(n)D(n)^2$.
We also need the classical upper bound of Vershik and Kerov
\cite{VK85}, restated in \cite[(2.4)]{PPY19}: there are universal
$v>0$ and $r_0$ such that
\begin{equation}\label{eq:VK-upper}
 \Delta(\rho)\ge v\sqrt r
       \qquad(\rho\vdash r,\ r\ge r_0).
\end{equation}

\begin{lemma}\label{lem:balance}
For every $A>0$ there is $B>1$, depending only on $A$, such that every
$\lambda\in\HH_A(n)$ is $B$-balanced.
\end{lemma}
\begin{proof}
Put $L:=\lambda_1$ and
$\bar\lambda:=(\lambda_2,\lambda_3,\ldots)\vdash n-L$.
Choosing the entries in the first row of a standard tableau shows that
\[
 f^\lambda\le\binom nL f^{\bar\lambda}
             \le\binom nL\sqrt{(n-L)!}.
\]
Consequently,
\begin{equation}\label{eq:row-removal}
 \frac{f^\lambda}{\sqrt{n!}}
       \le\frac{\sqrt{(n)_L}}{L!}
       \le\left(\frac{e\sqrt n}{L}\right)^L.
\end{equation}
Choose $B>e$ with $B\log(B/e)>A$. If $L>B\sqrt n$, the last
expression is smaller than $e^{-A\sqrt n}$, a contradiction.
Apply the same argument to $\lambda'$, whose dimension is unchanged.
\end{proof}

\subsection{Orthogonality and conjugacy classes}
Symmetric-group characters are real. The character formula for a
Kronecker coefficient therefore gives
\begin{equation}\label{eq:character-sum}
 R_n(\lambda,\mu,\nu)-1
       =\sum_{\pi\ts\in\ts S_n\setminus\{e\}}
                      r_\lambda(\pi)r_\mu(\pi)r_\nu(\pi).
\end{equation}
Row and column orthogonality
\cite[Theorems~1.9.3 and~1.10.3]{Sag01} imply
\begin{equation}\label{eq:orthogonality}
 \sum_{\pi\ts\in\ts S_n}|r_\lambda(\pi)|^2=e^{2\Delta(\lambda)},
 \qquad
 |r_\lambda(\pi)|\le
           \frac{e^{\Delta(\lambda)}}{\sqrt{|C_\pi|}},
\end{equation}
where $C_\pi$ is the conjugacy class of $\pi$.

\begin{lemma}\label{lem:class-size}
Let $s:=s(\pi)$. Then
\begin{equation}\label{eq:class-size}
 |C_\pi|\ge\left(\frac{\sqrt n}{e}\right)^s.
\end{equation}
\end{lemma}
\begin{proof}
If $m_j$ is the number of $j$-cycles for $j\ge2$, then
\[
 |C_\pi|=\frac{(n)_s}{\prod_{j\ts\ge\ts 2}j^{m_j}m_j!},
 \qquad
 \prod_{j\ts\ge\ts 2}j^{m_j}m_j!
       \le\prod_{j\ts\ge\ts 2}(jm_j)^{m_j}\le s^{s/2}.
\]
The geometric mean of the largest $s$ factors in $n!$ is at least
$(n!)^{1/n}\ge n/e$, so $(n)_s\ge(n/e)^s$. Since $s\le n$, the
claim follows. The case $s=0$ is immediate.
\end{proof}

\subsection{Content specializations}
We record three probability weights that will turn a dimension bound
into content estimates. Schur--Weyl duality
\cite[Chapter~I, Appendix~A, \S\S5, 8]{Mac95} gives, for positive integral
$d$,
\begin{equation}\label{eq:SW}
 Q_{n,d}(\lambda):=\frac{f^\lambda s_\lambda(1^d)}{d^n}
       =\Pl_n(\lambda)\prod_{b\ts\in\ts \lambda}\left(1+\frac{c(b)}d\right),
 \qquad \sum_{\lambda\ts\vdash\ts n}Q_{n,d}(\lambda)=1.
\end{equation}
The ordinary and dual Cauchy identities
\cite[Chapter~I, (4.3) and (4.3$'$)]{Mac95}, under the specializations in
\eqref{eq:hooks}, give
\begin{equation}\label{eq:Cauchy}
 \sum_{\lambda\ts\vdash\ts n}s_\lambda(1^z)s_\lambda(1^w)
      =\frac{(zw)^{\overline n}}{n!},
 \qquad
 \sum_{\lambda\ts\vdash\ts n}s_\lambda(1^D)s_{\lambda'}(1^D)
      =\binom{D^2}{n},
\end{equation}
where $x^{\overline n}:=x(x+1)\cdots(x+n-1)$ and $D$ is a positive
integer. These identities follow by taking the coefficient of $t^n$
in $(1-t)^{-zw}$ and $(1+t)^{D^2}$, respectively.
Taking $z=iD$, $w=-iD$ in the first identity gives the probability weight
\begin{equation}\label{eq:Qplus}
 Q^+_{n,D}(\lambda)
  :=\frac{n!\,|s_\lambda(1^{iD})|^2}{(D^2)^{\overline n}}
   =\Pl_n(\lambda)
     \frac{\prod_{b\ts\in\ts \lambda}(1+c(b)^2/D^2)}
          {\prod_{j\ts=\ts 0}^{n-1}(1+j/D^2)}.
\end{equation}
When $n\le D^2$, the second identity gives another probability weight,
\begin{equation}\label{eq:Qminus}
 Q^-_{n,D}(\lambda)
  :=\frac{n!\,s_\lambda(1^D)s_{\lambda'}(1^D)}{(D^2)_n}
   =\Pl_n(\lambda)
     \frac{\prod_{b\ts\in\ts \lambda}(1-c(b)^2/D^2)}
          {\prod_{j\ts=\ts 0}^{n-1}(1-j/D^2)}.
\end{equation}
Its first expression is nonnegative and vanishes outside the $D$ by $D$
square. Thus nonnegativity does not require interpreting the signs of
individual factors outside that square.

The purpose of \eqref{eq:Qplus} and \eqref{eq:Qminus} is to obtain
upper and lower bounds for the second content moment. For a
$B$-balanced diagram and $D/\sqrt n\to\infty$, expansion of their
logarithms gives
\[
 \log\frac{Q^\pm_{n,D}(\lambda)}{\Pl_n(\lambda)}
 =\pm\frac{\sum_{b\ts\in\ts \lambda}c(b)^2-\binom n2}{D^2}
       +O_B(n^3/D^4).
\]
The denominators in the two weights supply the centering term
$\binom n2$. Since each weight is at most one, a lower bound on
$\Pl_n(\lambda)$ bounds the displayed logarithm from above.
The plus sign therefore controls the positive deviation of the second
content moment, and the minus sign controls its negative deviation.
We use both bounds in Lemma~\ref{lem:contents}.

\subsection{The Stanley--F\'eray formula}\label{subsec:SF}
For $\lambda\vdash n$ and $\alpha,\beta\in S_s$, write $C(\alpha)$
and $C(\beta)$ for their sets of cycles, including one-cycles. A
\defn{coloring} is a pair of maps
\[
 u:C(\alpha)\longrightarrow\{1,\ldots,\lambda_1\},
 \qquad
 v:C(\beta)\longrightarrow\{1,\ldots,\ell(\lambda)\}.
\]
Thus $u$ assigns a column number to each $\alpha$-cycle, and $v$
assigns a row number to each $\beta$-cycle. The coloring is
\defn{compatible with $\lambda$} if
\[
 A\cap B\ne\varnothing
       \quad\Longrightarrow\quad
       (v(B),u(A))\in\lambda,
 \qquad A\in C(\alpha),\ B\in C(\beta).
\]
Here cycles are viewed as sets of labels when taking their
intersection. Equivalently, the condition is $u(A)\le\lambda_{v(B)}$
whenever the two cycles intersect. Define $N^\lambda(\alpha,\beta)$
to be the number of compatible pairs $(u,v)$. Neither map is required
to be injective: different cycles may receive the same row or column
number. For example, when $\alpha=\beta=\mathrm{id}_s$, one may choose
a cell independently for each label, so
$N^\lambda(\mathrm{id}_s,\mathrm{id}_s)=n^s$.

It is useful to encode these conditions in a \defn{bipartite cycle
graph} $G(\alpha,\beta)$. Its two vertex sets are disjoint copies of
$C(\alpha)$ and $C(\beta)$. For each label $a\in\{1,\ldots,s\}$,
draw an edge labeled $a$ between the $\alpha$-cycle and the
$\beta$-cycle containing $a$. Thus parallel edges are allowed.
Coloring the vertices means assigning column numbers to the
$\alpha$-vertices and row numbers to the $\beta$-vertices. Each edge
then specifies a cell, whose row and column are the colors of its two
endpoints; compatibility requires this cell to belong to $\lambda$.
These are numerical assignments, with no requirement that adjacent
vertices have distinct colors. The connected components of this graph
correspond to the orbits of $\langle\alpha,\beta\rangle$: moving
between labels in a common cycle generates exactly the same
equivalence relation as applying $\alpha$ and $\beta$.

\begin{example}\label{ex:coloring}
Take $\lambda:=(3,2)$, $\alpha:=(1\,2)(3)$, and
$\beta:=(1)(2\,3)$. Write $A_{12},A_3$ for the two
$\alpha$-vertices and $B_1,B_{23}$ for the two $\beta$-vertices.
The graph is the path
\[
 B_1\ \overset{1}{\longleftrightarrow}\ A_{12}
 \ \overset{2}{\longleftrightarrow}\ B_{23}
 \ \overset{3}{\longleftrightarrow}\ A_3.
\]
Write $x:=u(A_{12})$, $y:=u(A_3)$, $r:=v(B_1)$, and $w:=v(B_{23})$.
The three edge conditions are
\[
 (r,x),\ (w,x),\ (w,y)\in\lambda.
\]
For fixed rows $r,w\in\{1,2\}$, there are
$\min(\lambda_r,\lambda_w)$ choices for $x$ and $\lambda_w$
choices for $y$. Consequently,
\[
\begin{array}{c|cccc}
 (r,w)&(1,1)&(1,2)&(2,1)&(2,2)\\ \hline
 \#\{(x,y)\}&3\cdot3&2\cdot2&2\cdot3&2\cdot2
\end{array}
 \qquad
 N^{(3,2)}(\alpha,\beta)=9+4+6+4=23.
\]
For instance, $r=w=1$ and $x=y=2$ is allowed: all three edges
specify the same cell $(1,2)$. In contrast, $w=2$ and $x=3$ is
forbidden because $(2,3)\notin\lambda$.
\end{example}

The character formula in \cite[Theorem~2]{FS11} is
\begin{equation}\label{eq:SF}
 \Sigma_\lambda(\pi)
     =\sum_{\substack{\alpha,\beta\ts\in\ts S_s\\\alpha\beta\ts=\ts \pi}}
                   (-1)^{|\alpha|}N^\lambda(\alpha,\beta),
       \qquad s=s(\pi).
\end{equation}
Let $O_1,\ldots,O_b$ be the orbits of $\langle\alpha,\beta\rangle$
on the labels, and put
$\alpha_h:=\alpha|_{O_h}$ and $\beta_h:=\beta|_{O_h}$, relabeling
each orbit when needed. Every cycle lies in one orbit, and cycles
in different orbits do not intersect. Thus a compatible coloring
is exactly an independent choice of a compatible coloring on each
orbit; colors may be reused in different orbits. Consequently,
\[
 N^\lambda(\alpha,\beta)
     =\prod_{h\ts=\ts 1}^b N^\lambda(\alpha_h,\beta_h),
 \qquad
 (-1)^{|\alpha|}=\prod_{h\ts=\ts 1}^b(-1)^{|\alpha_h|}.
\]
The sign identity follows because transposition length is additive
on disjoint label sets. This is what we mean by
\defn{multiplicative over the orbits}. It lets us sum signed
contributions on selected orbits before taking absolute values.

\subsection{Restriction and induction norms}
For $\rho\vdash r$ and $0\le a\le r$, define
\begin{equation}\label{eq:Bbeta}
 B_\rho(a):=\sum_{\substack{\alpha\ts\vdash\ts r-a\\\beta\ts\vdash\ts a}}
                         (c^\rho_{\alpha\beta})^2,
 \qquad \beta_\rho(a):=\frac{B_\rho(a)}{\binom ra}.
\end{equation}
The squared character norm of a restriction gives
\begin{equation}\label{eq:restriction-probability}
 \beta_\rho(a)
       =\frac1{r!}\sum_{\pi\ts\in\ts S_{r-a}\times S_a}|\chi^\rho(\pi)|^2
       =Q_\rho(S_{r-a}\times S_a)\le1,
 \qquad Q_\rho(\pi):=\frac{|\chi^\rho(\pi)|^2}{r!}.
\end{equation}
Here $Q_\rho$ is a central probability measure by orthogonality.

\begin{prop}[An induction-norm bound]\label{prop:mackey}
If $\mu\vdash k$, $\nu\vdash l$, and $n=k+l$, then
\begin{equation}\label{eq:mackey-bound}
 (c^\lambda_{\mu\nu})^2
       \le\sum_{\tau\ts\vdash\ts n}(c^\tau_{\mu\nu})^2
       \le\sum_{a\ts=\ts 0}^{\min(k,l)}B_\mu(a)B_\nu(a)
       \le\binom nk.
\end{equation}
\end{prop}
\begin{proof}
We use \cite[Lemma~5.1]{PPY19}, writing $\tau$ for its outer
summation index and $\xi$ for its $\delta$. In its four-partition sum, a
nonzero term forces $|\gamma|=|\xi|$ by the size conditions for
$c^\mu_{\alpha\gamma}$ and $c^\mu_{\alpha\xi}$. Writing
$a:=|\gamma|=|\xi|$, the other sizes are $|\alpha|=k-a$ and
$|\beta|=l-a$, so $0\le a\le\min(k,l)$. Renaming their summation
indices and grouping by $a$ gives
\begin{equation}\label{eq:mackey-identity}
 \sum_{\tau\ts\vdash\ts n}(c^\tau_{\mu\nu})^2
  =\sum_{a\ts=\ts 0}^{\min(k,l)}\ \sum_{\substack{\alpha\ts\vdash\ts k-a,\ \beta\ts\vdash\ts l-a\\
                           \gamma,\xi\ts\vdash\ts a}}
       c^\mu_{\alpha\gamma}c^\mu_{\alpha\xi}
       c^\nu_{\beta\gamma}c^\nu_{\beta\xi}.
\end{equation}
For fixed $a$, let
$M_\mu:=(c^\mu_{\alpha\gamma})_{\alpha,\gamma}$ and
$M_\nu:=(c^\nu_{\beta\gamma})_{\beta,\gamma}$. Their common
column index runs over partitions of $a$. The inner sum is
$\tr(XY)$ for $X:=M_\mu^{\T}M_\mu$ and
$Y:=M_\nu^{\T}M_\nu$. These matrices are positive semidefinite,
and their traces are $B_\mu(a)$ and $B_\nu(a)$. Therefore
$\tr(XY)\le\tr(X)\tr(Y)$ gives the middle bound in
\eqref{eq:mackey-bound}; its first bound keeps one nonnegative
summand. Finally, \defn{Vandermonde's identity} states

\begin{equation}\label{eq:vandermonde}
 \sum_{a\ts=\ts 0}^{\min(k,l)}\binom ka\binom la=\binom{k+l}{k}.
\end{equation}
Together with $B_\rho(a)\le\binom{|\rho|}{a}$ from
\eqref{eq:restriction-probability}, this gives
\[
 \sum_a B_\mu(a)B_\nu(a)
       \le\sum_a\binom ka\binom la=\binom nk.
\]
\end{proof}

By Frobenius reciprocity, restriction has the decomposition
\[
 \Res_{S_k\times S_l}^{S_n}V^\lambda
   \cong\bigoplus_{\substack{\alpha\ts\vdash\ts k\\\beta\ts\vdash\ts l}}
           (V^\alpha\boxtimes V^\beta)^{\oplus c^\lambda_{\alpha\beta}}.
\]
Taking dimensions and keeping only the summand indexed by
$(\mu,\nu)$ gives the first inequality below. Normalizing it by
$\sqrt{n!}$ and taking logarithms gives the second:
\begin{equation}\label{eq:dimension-restriction}
 f^\lambda\ge c^\lambda_{\mu\nu}f^\mu f^\nu,
 \qquad
 \Delta(\lambda)\le\Delta(\mu)+\Delta(\nu)
                             +\delta(\lambda;\mu,\nu)
       \quad(c^\lambda_{\mu\nu}>0).
\end{equation}
More explicitly, the identity
$\sqrt{n!}=\sqrt{\binom nk}\sqrt{k!}\sqrt{l!}$ gives, when
$c^\lambda_{\mu\nu}>0$,
\[
 e^{-\Delta(\lambda)}
   \ge\frac{c^\lambda_{\mu\nu}}{\sqrt{\binom nk}}
          \frac{f^\mu}{\sqrt{k!}}\frac{f^\nu}{\sqrt{l!}}
   =e^{-\delta(\lambda;\mu,\nu)-\Delta(\mu)-\Delta(\nu)}.
\]
Taking negative logarithms proves the asserted deficit inequality.

\section{Uniform Kronecker asymptotics}\label{sec:kronecker}

We prove Theorem~\ref{thm:uniform} by showing that the cubic absolute
sum of each normalized character tends to zero, uniformly in the
large-dimension window. First, content identities give cancellation
at transpositions and $3$-cycles. We retain this cancellation in a
factorization bound, sum over permutations of small support, and use
orthogonality for the remaining support sizes. H\"older's inequality
then gives the relative asymptotic and its corollary.

Recall that $\HH_A(n):=\{\lambda\vdash n:
f^\lambda\ge\sqrt{n!}e^{-A\sqrt n}\}$ and that
$r_\lambda(\pi):=\chi^\lambda(\pi)/f^\lambda$. Also recall
$s(\pi):=|\supp\pi|$ and
\[
 \Sigma_\lambda(\pi):=(n)_{s(\pi)}r_\lambda(\pi),
 \qquad
 \Sigma_\lambda(j):=(n)_j r_\lambda((1\,2\,\cdots\,j))
      \quad(2\le j\le n),
\]
where $(n)_j:=n(n-1)\cdots(n-j+1)$. Thus
$\Sigma_\lambda(j)$ denotes the normalized character on one
$j$-cycle, with the other labels fixed. We use the probability weights $Q_{n,d}$ and $Q^\pm_{n,D}$ defined in
\eqref{eq:SW}, \eqref{eq:Qplus} and \eqref{eq:Qminus}.

\subsection{Cancellation in the first two content moments}
Put $M_j(\lambda):=\sum_{b\ts\in\ts \lambda}c(b)^j$.
\begin{lemma}\label{lem:contents}
For fixed $A>0$, uniformly over $\lambda\in\HH_A(n)$,
\begin{equation}\label{eq:content-estimates}
 M_1(\lambda)=O_A(n^{5/4}),
 \qquad M_2(\lambda)-\binom n2=O_A(n^{7/4}).
\end{equation}
Consequently,
\begin{equation}\label{eq:small-cycle-estimates}
 \Sigma_\lambda(2)=O_A(n^{5/4}),
 \qquad \Sigma_\lambda(3)=O_A(n^{7/4}).
\end{equation}
\end{lemma}
\begin{proof}
By Lemma~\ref{lem:balance}, $|c(b)|\le B\sqrt n$ for $B:=B(A)$.
Also $\log\Pl_n(\lambda)\ge-2A\sqrt n$.
In \eqref{eq:SW}, take $d:=\lceil n^{3/4}\rceil$. For large $n$, all
$|c(b)/d|\le1/2$. Since $Q_{n,d}(\lambda)\le1$ and
$\log(1+x)\ge x-x^2$ for $|x|\le1/2$,
\[
 \frac{M_1(\lambda)}d-\frac{M_2(\lambda)}{d^2}
       \le-\log\Pl_n(\lambda)\le2A\sqrt n.
\]
The crude bound $M_2(\lambda)\le B^2n^2$ gives
$M_1(\lambda)\le O_A(n^{5/4})$. Replacing $\lambda$ by its conjugate
$\lambda'$ leaves the dimension and $M_2$ unchanged, while replacing
$M_1$ by $-M_1$. Applying the same upper bound to $\lambda'$ therefore
gives the matching lower bound for $M_1(\lambda)$.

Next take $D:=\lceil n^{5/8}\rceil$ in \eqref{eq:Qplus} and
\eqref{eq:Qminus}. For large $n$, the diagram is contained in the $D$
by $D$ square and all displayed factors are positive. Using
$M_4(\lambda)\le B^4n^3$ and $\sum_{j\ts<\ts n}j^2=O(n^3)$, the Taylor
expansions of the logarithms give
\[
 \log Q^\pm_{n,D}(\lambda)
   =\log\Pl_n(\lambda)
      \pm\frac{M_2(\lambda)-\binom n2}{D^2}
      +O_B(n^3/D^4).
\]
Both left sides are nonpositive. Combining the two signs yields
\[
 \left|M_2(\lambda)-\binom n2\right|
       \le2AD^2\sqrt n+O_B(n^3/D^2)=O_A(n^{7/4}).
\]

To pass to characters, let $J_i:=\sum_{j\ts<\ts i}(j\,i)$ be the
Jucys--Murphy elements, with $J_1:=0$. In the seminormal basis
$\{v_T\}$ of $V^\lambda$, indexed by standard Young tableaux $T$ of
shape $\lambda$, their eigenvalue formula is
\[
 J_i v_T=c_T(i)v_T,
\]
where $c_T(i)$ is the content of the cell containing $i$; see
\cite[\S3, (3.8)]{Mur81}. Since each cell contains one label, it follows that
\[
 \left(\sum_iJ_i\right)v_T=M_1(\lambda)v_T,
 \qquad
 \left(\sum_iJ_i^2\right)v_T=M_2(\lambda)v_T.
\]
Let $\mathcal C_2$ and $\mathcal C_3$ be the conjugacy classes of transpositions and
$3$-cycles, and put $K_j:=\sum_{\pi\ts\in\ts \mathcal C_j}\pi$. In the group
algebra one has
\[
 \sum_iJ_i=K_2,
 \qquad
 \sum_iJ_i^2=\binom n2e+K_3.
\]
The first identity counts each transposition once. For the second,
expand
\[
 J_i^2=(i-1)e+
       \sum_{1\ts\le\ts j\ts<\ts k\ts<\ts i}
          \bigl((j\,i)(k\,i)+(k\,i)(j\,i)\bigr).
\]
For each triple of labels, its largest label determines $i$, and the
two products give its two $3$-cycles. Thus every $3$-cycle occurs
once, while the identity has coefficient
$\sum_i(i-1)=\binom n2$.

For any conjugacy class $C$, its class sum is central and acts on
$V^\lambda$ by a scalar. Taking its trace shows that the scalar is
$|C|\chi^\lambda(C)/f^\lambda=|C|r_\lambda(C)$. Comparing this with
the content eigenvalues above gives
\[
 |\mathcal C_2|r_\lambda(2)=M_1(\lambda),
 \qquad
 |\mathcal C_3|r_\lambda(3)=M_2(\lambda)-\binom n2.
\]
Finally, $|\mathcal C_2|=(n)_2/2$ and $|\mathcal C_3|=(n)_3/3$. Multiplying by
$(n)_2$ and $(n)_3$, respectively, and using the definition of
$\Sigma_\lambda$ yields
\begin{equation}\label{eq:content-character}
 \Sigma_\lambda(2)=2M_1(\lambda),
 \qquad
 \Sigma_\lambda(3)=3\left(M_2(\lambda)-\binom n2\right).
\end{equation}
This proves \eqref{eq:small-cycle-estimates}.
\end{proof}

\subsection{A character bound retaining transposition cancellation}
The usual balanced-diagram character bound does not by itself give a
summable cubic estimate at transpositions. We need to preserve their
signed contribution in \eqref{eq:SF}.
Here $B$-balanced means $\lambda_1,\ell(\lambda)\le B\sqrt n$.
Recall also $s(\pi):=|\supp\pi|$ and
$\Sigma_\lambda(\pi):=(n)_{s(\pi)}r_\lambda(\pi)$, where
$(n)_s:=n(n-1)\cdots(n-s+1)$; the coloring count
$N^\lambda(\alpha,\beta)$ is defined in Section~\ref{subsec:SF}.

\begin{lemma}\label{lem:refined}
For every $B\ge1$ there is $K>0$, depending only on $B$, with the following property.
Suppose $\lambda\vdash n$, with $n\ge3$, is $B$-balanced, and $\pi$ has support
$s\le\sqrt n$ and $q$ nontrivial cycles, of which $m$ are
transpositions. Put $a_\lambda:=|\Sigma_\lambda(2)|/n^{3/2}$. Then
\begin{equation}\label{eq:refined}
 |r_\lambda(\pi)|\le
       K^s n^{-(s-q)/2}
          \left(a_\lambda+\frac{s}{\sqrt n}\right)^m.
\end{equation}
\end{lemma}
\begin{proof}
If $s=0$, the assertion is immediate. We therefore assume $s\ge2$.

\smallskip
\noindent\emph{The coloring bound.}
We first prove an auxiliary estimate for a permutation $\tau$ acting
on $t$ labels, all of which it moves. After relabeling its support as
$\{1,\ldots,t\}$, we regard $\tau$ as an element of $S_t$.
Later, $\tau$ will be obtained from $\pi$ by removing some isolated
transpositions. In this auxiliary argument, $q$ denotes the number
of cycles of $\tau$. For a factorization $\alpha\beta=\tau$ in
$S_t$, set
\[
 b:=\#\operatorname{Orb}\langle\alpha,\beta\rangle,
 \qquad V:=|C(\alpha)|+|C(\beta)|,
 \qquad o:=t-V+b,
\]
where all cycle counts include one-cycles. Let $L$ satisfy
$\lambda_1,\ell(\lambda)\le L$ and $t\le L\le n$.
Apply the first inequality of \cite[Proposition~12, (18)]{FS11},
with their parameters $l,A,\pi$ equal to our $t,L,\tau$ and their
factors equal to $\alpha,\beta$. All its hypotheses hold by the
choice of $L$. It gives
\begin{equation}\label{eq:color-bound}
 N^\lambda(\alpha,\beta)
   \le n^bL^{V-2b}
   =L^{t-q}n^qL^{-o}(L/n)^{q-b}.
\end{equation}
The equality uses $V=t+b-o$. We retain the orbit dependence in
the first bound of that proposition: $(L/n)^{q-b}$ is the gain
when several cycles of $\tau$ belong to one generated orbit.

\smallskip
\noindent\emph{Reducing and counting factorizations.}
The reduction in \cite[Lemma~13 and its proof]{FS11} replaces
$(\alpha,\beta)$ by $(\alpha',\beta')$ with the same product and
orbit partition, such that every $\alpha'$-cycle is contained in
a $\beta'$-cycle. Fix one such reduced pair for each input pair.
The $\beta'$-cycles are then precisely the $b$ generated orbits,
and the lemma gives
\begin{equation}\label{eq:reduction}
 \alpha'\beta'=\tau,
 \qquad |\alpha'|=o,
 \qquad |C(\beta')|=b,
 \qquad |\beta'|=|\beta'\beta^{-1}|+|\beta|.
\end{equation}
In particular, $|\beta'|=t-b$, by the transposition-length formula
of Section~\ref{subsec:characters}. The last equality says that
$\beta$ lies on a shortest transposition path from the identity to
$\beta'$: a sequence of left multiplications by transpositions
using the minimum possible number of steps.

For a fixed $\beta'$ with cycle lengths $d_1,\ldots,d_b$,
\cite[\S5.1]{FS11} identifies the possible $\beta$ with noncrossing
partitions of its cycles. On each cycle, arrange the labels in
cyclic order on a circle; the blocks are the $\beta$-cycles, their
chords do not cross, and each block inherits the cyclic order.
The choices on distinct cycles are independent, giving at most
\[
 \prod_{h\ts=\ts 1}^b C_{d_h}\le4^{d_1+\cdots+d_b}=4^t,
 \qquad C_d:=\frac1{d+1}\binom{2d}{d}\le4^d.
\]
Here $C_d$ is the $d$th Catalan number. Each $\beta$ determines
$\alpha=\tau\beta^{-1}$, so a reduced pair has at most $4^t$
preimages. Since $\alpha'$ determines
$\beta'=(\alpha')^{-1}\tau$, the number of reduced pairs with
$|\alpha'|=j$ is bounded by \cite[Lemma~14]{FS11}:
\begin{equation}\label{eq:length-count}
 \#\{\sigma\in S_t:|\sigma|=j\}\le t^{2j}/j!.
\end{equation}
Finally, $\tau$ has $q$ cycles, so the triangle inequality gives
\begin{equation}\label{eq:orbit-defect}
 t-q=|\tau|\ge|\beta'|-|\alpha'|=t-b-o,
 \qquad o\ge q-b.
\end{equation}

\smallskip
\noindent\emph{The cost of nonisolated transpositions.}
Fix a set $U$ of $u$ transposition cycles of $\tau$, and let
$\mathcal F_U(\tau)$ be the set of factorizations
$\alpha\beta=\tau$ for which none of the two-element supports of
these cycles is an orbit of $\langle\alpha,\beta\rangle$.
Every generated orbit is a union of cycles of $\tau$, because
$\tau=\alpha\beta$ preserves that orbit. If the $b$ orbits contain
$h_1,\ldots,h_b$ cycles of $\tau$, then
\[
 \sum_{a\ts=\ts 1}^b h_a=q,
 \qquad q-b=\sum_{a\ts=\ts 1}^b(h_a-1).
\]
An orbit containing a specified transposition must contain at least
one further $\tau$-cycle, so its $h_a$ is at least two. It contains
at most $h_a\le2(h_a-1)$ specified transpositions. Summing over
these orbits gives $u\le2(q-b)$, and therefore
\[
 q-b\ge v:=\lceil u/2\rceil,
 \qquad o\ge v
\]
by \eqref{eq:orbit-defect}. For example, two specified transpositions
can form one four-label orbit; this contributes one to $q-b$.
This explains the factor of two in the bound.

We now sum the absolute values of the individual terms in the
Stanley--F\'eray formula, restricted to $\mathcal F_U(\tau)$.
Since $N^\lambda(\alpha,\beta)$ is nonnegative, each such absolute
value is simply $N^\lambda(\alpha,\beta)$.
Group the factorizations by the length $j=|\alpha'|=o$ of their
reduced image. There are at most $t^{2j}/j!$ possible images and at
most $4^t$ preimages per image. In \eqref{eq:color-bound}, use
$L^{-o}=L^{-j}$ and $(L/n)^{q-b}\le(L/n)^v$, since $L/n\le1$.
It follows that
\begin{align}
 \sum_{(\alpha,\beta)\ts\in\ts \mathcal F_U(\tau)}
     \bigl|(-1)^{|\alpha|}N^\lambda(\alpha,\beta)\bigr|
 &=\sum_{(\alpha,\beta)\ts\in\ts \mathcal F_U(\tau)}
          N^\lambda(\alpha,\beta)\notag\\
 &\le4^tL^{t-q}n^q(L/n)^v
     \sum_{j\ts\ge\ts v}\frac{(t^2/L)^j}{j!}
       \notag\\
 &\le4^tL^{t-q}n^q e^{t^2/L}\frac{(t^2/n)^v}{v!}\notag\\
 &\le4^tL^{t-q}n^q e^{t^2/L}(t/\sqrt n)^u,
 \label{eq:residual-bound}
\end{align}
provided $t\le\sqrt n$.
For the second inequality, put $x:=t^2/L$ and use
\[
 \sum_{j\ts\ge\ts v}\frac{x^j}{j!}
 =x^v\sum_{h\ts\ge\ts 0}\frac{x^h}{(v+h)!}
 \le\frac{x^v}{v!}e^x,
\]
which follows from $(v+h)!\ge v!h!$; multiplication by $(L/n)^v$
then gives $(t^2/n)^v/v!$. For the last inequality in
\eqref{eq:residual-bound}, use $2v\ge u$, $t/\sqrt n\le1$, and
$1/v!\le1$.

\smallskip
\noindent\emph{Keeping the isolated transpositions signed.}
Return to $\pi$ and the parameters $s,q,m$ in the statement. Let
$T$ be its set of $m$ transposition cycles. In \eqref{eq:SF}, group
the factorizations according to the exact subset $J\subseteq T$
whose two-element supports are generated orbits. On each such
isolated support $\{i,j\}$, the only factorizations are
$(\mathrm{id},(i\,j))$ and $((i\,j),\mathrm{id})$. Their signed
sum is $\Sigma_\lambda(2)$ by \eqref{eq:SF} for two labels. In
particular, we sum these two terms before taking absolute values.
The coloring count is multiplicative over components: colors in
different components are unconstrained and may coincide. The signs
are also multiplicative, since transposition length is additive on
disjoint supports.

Write $|J|=m-u$. Removing its isolated supports leaves a
permutation $\tau$ with
\[
 t:=s-2(m-u),\qquad q':=q-(m-u).
\]
It has $u$ specified transposition cycles, namely $U:=T\setminus J$,
none of which is an individual orbit in the residual factorization.
Thus the total contribution of the group indexed by $J$ is exactly
\[
 \Sigma_\lambda(2)^{m-u}
 \sum_{(\alpha_0,\beta_0)\ts\in\ts \mathcal F_U(\tau)}
      (-1)^{|\alpha_0|}N^\lambda(\alpha_0,\beta_0).
\]
For $t>0$, apply \eqref{eq:residual-bound} with $q'$ in place of
$q$ and $L:=\min(B\sqrt n,n)$. This $L$ bounds both diagram
dimensions, satisfies $t\le L\le n$, and obeys
\[
 e^{t^2/L}\le e^t,
 \qquad L^{t-q'}\le B^{t-q'}n^{(t-q')/2}
                         \le B^t n^{(t-q')/2}.
\]
The absolute value of the residual sum is consequently at most
\[
 (4eB)^t n^{(t+q')/2}(t/\sqrt n)^u.
\]
When $t=0$, the residual factor is $\Sigma_\lambda(\varnothing)=1$
and $u=0$; we treat this empty case directly.
For a nonempty residual permutation, multiply the bound by
$|\Sigma_\lambda(2)|^{m-u}
 =a_\lambda^{m-u}n^{3(m-u)/2}$. The powers of $n$ combine because
\[
 \frac{t+q'}2+\frac32(m-u)=\frac{s+q}2,
\]
as follows by substituting the formulas for $t$ and $q'$.
Also $(4eB)^t\le(4eB)^s$ and $t^u\le s^u$.
There are $\binom m{m-u}=\binom mu$ possible subsets $J$ of this
size. Including the empty residual case and summing over $u$ gives
\[
 |\Sigma_\lambda(\pi)|
   \le(4eB)^s n^{(s+q)/2}
       \sum_{u\ts=\ts 0}^m\binom mu a_\lambda^{m-u}(s/\sqrt n)^u.
\]
The binomial theorem identifies the sum with
$(a_\lambda+s/\sqrt n)^m$. Finally, the first $s$ factors of $n!$
have geometric mean at least that of all $n$ factors, and
$n!\ge(n/e)^n$, so $(n)_s\ge(n/e)^s$. Dividing by $(n)_s$ in
the definition of $\Sigma_\lambda$ gives
\[
 |r_\lambda(\pi)|\le
 (4e^2B)^s n^{-(s-q)/2}
       (a_\lambda+s/\sqrt n)^m.
\]
This is \eqref{eq:refined} with $K:=4e^2B$.
\end{proof}

\begin{lemma}\label{lem:three-cycles}
For fixed $B\ge1$ and $q\ge1$, uniformly over $B$-balanced
$\lambda\vdash n$ with $n\ge3q$,
\begin{equation}\label{eq:three-cycles}
 \Sigma_\lambda(3^q)=\Sigma_\lambda(3)^q+O_{B,q}(n^{2q-1}).
\end{equation}
\end{lemma}
\begin{proof}
Relabel the support as $\{1,\ldots,3q\}$, and write
$\pi=c_1\cdots c_q$ for its disjoint $3$-cycles, with supports
$I_1,\ldots,I_q$. In \eqref{eq:SF}, the factors lie in $S_{3q}$.
For a factorization $\alpha\beta=\pi$, every
orbit of $\langle\alpha,\beta\rangle$ is a union of these supports:
the group contains $\pi$, which acts transitively on each $I_h$.
Consequently, the number $b$ of generated orbits is at most $q$,
and equality holds precisely when the orbits are $I_1,\ldots,I_q$.
Every other factorization joins at least two of the $3$-cycles in
one orbit, so it has $b\le q-1$.

When the orbits are exactly $I_1,\ldots,I_q$, both factors restrict
to independent factorizations $\alpha_h\beta_h=c_h$ on the three
labels of $I_h$. Conversely, any such choices give these exact
orbits, since each restricted generated group contains the
transitive cycle $c_h$. The sign and coloring count multiply over
the supports, as explained in Section~\ref{subsec:SF}; in particular,
colors may be reused on different supports. Thus this part of
\eqref{eq:SF} is
\[
 \prod_{h\ts=\ts 1}^q
   \left(\sum_{\alpha_h\beta_h\ts=\ts c_h}
            (-1)^{|\alpha_h|}N^\lambda(\alpha_h,\beta_h)\right)
   =\prod_{h\ts=\ts 1}^q\Sigma_\lambda(c_h)
   =\Sigma_\lambda(3)^q,
\]
by applying \eqref{eq:SF} on each three-element support.

For the remaining factorizations, use
$V:=|C(\alpha)|+|C(\beta)|$ and $o:=3q-V+b$, as in the proof of
Lemma~\ref{lem:refined}. Equation~\eqref{eq:orbit-defect} gives
$o\ge q-b$, and hence
\[
 |C(\alpha)|+|C(\beta)|=3q+b-o\le2q+2b\le4q-2.
\]
Each vertex has at most $B\sqrt n$ possible colors, so each
remaining coloring count is at most
$(B\sqrt n)^{4q-2}=O_{B,q}(n^{2q-1})$. There are at most $(3q)!$
factorizations in total, because choosing $\beta$ determines
$\alpha=\pi\beta^{-1}$. Summing their absolute values proves the
error estimate.
\end{proof}

\subsection{The cubic character sum}
\begin{prop}\label{prop:cubic}
For every fixed $A>0$,
\begin{equation}\label{eq:cubic}
 \max_{\lambda\ts\in\ts \HH_A(n)}
       \sum_{\pi\ts\in\ts S_n\setminus\{e\}}|r_\lambda(\pi)|^3
       \longrightarrow0.
\end{equation}
\end{prop}
\begin{proof}
Fix $B:=B(A)$ from Lemma~\ref{lem:balance} and $K:=K(B)$ from
Lemma~\ref{lem:refined}. Put
\[
 a_\lambda:=|\Sigma_\lambda(2)|/n^{3/2},
 \qquad \eta_\lambda:=n^{1/2}a_\lambda^3.
\]
By Lemma~\ref{lem:contents},
$a_\lambda=O_A(n^{-1/4})$ and $\eta_\lambda=O_A(n^{-1/4})$.

Define the contribution from support exactly $s$ by
\[
 T_{\lambda,s}:=
   \sum_{\substack{\pi\ts\in\ts S_n\\s(\pi)\ts=\ts s}}|r_\lambda(\pi)|^3.
\]
For multiplicities $\mathbf m:=(m_j)_{j\ts\ge\ts 2}$ satisfying
$\sum_{j\ts\ge\ts 2}jm_j=s$, let $\mathcal C_{\mathbf m}$ denote the
conjugacy class of type $1^{n-s}2^{m_2}3^{m_3}\cdots$. Its size is
$(n)_s/\prod_{j\ts\ge\ts 2}j^{m_j}m_j!$.
For $s\le\sqrt n$, Lemma~\ref{lem:refined} gives the following
bound on the sum over this entire class:
\begin{equation}\label{eq:type-majorant}
 \begin{aligned}
 \sum_{\pi\ts\in\ts \mathcal C_{\mathbf m}}|r_\lambda(\pi)|^3
   &=\frac{(n)_s}{\prod_{j\ts\ge\ts 2}j^{m_j}m_j!}
           |r_\lambda(2^{m_2}3^{m_3}\cdots)|^3\\
   &\le\prod_{j\ts\ge\ts 2}\frac{b_j(n,s)^{m_j}}{m_j!},
 \end{aligned}
\end{equation}
where
\[
 b_2(n,s):=\frac{K^6}{2}n^{1/2}(a_\lambda+s/\sqrt n)^3,
 \qquad
 b_j(n,s):=\frac{K^{3j}}j n^{(3-j)/2}\quad(j\ge3).
\]
Indeed, the exponent of $n$ after multiplying by the class size is
$(3\sum_jm_j-s)/2$.

For a fixed $s$, summing the right side of
\eqref{eq:type-majorant} over types gives
\[
 T_{\lambda,s}\le[z^s]F_{\lambda,n,s}(z),
 \qquad
 F_{\lambda,n,s}(z):=
          \exp\left(\sum_{j\ts\ge\ts 2}b_j(n,s)z^j\right).
\]
Here $[z^s]$ denotes coefficient extraction; expansion of each
$\exp(b_jz^j)$ selects $m_j$ cycles of length $j$.
Although $b_2$ depends on $s$, we use this coefficient bound separately
for each $s$. The inequality $(x+y)^3\le4(x^3+y^3)$ gives
\[
 b_2(n,s)\le2K^6(\eta_\lambda+s^3/n).
\]
Choose $\varepsilon>0$ so small that
$8K^6\varepsilon^2\le\tfrac12\log2$. For
$s\le\varepsilon\sqrt n$ and large $n$,
\[
 4b_2(n,s)\le8K^6+\tfrac12s\log2,
 \qquad
 \sum_{j\ts\ge\ts 3}2^jb_j(n,s)
      =\sum_{j\ts\ge\ts 3}\frac{(2K^3)^j}{j}n^{(3-j)/2}=O_B(1).
\]
The series converges at $z=2$ for $n>(2K^3)^2$. We evaluate
the generating function $F_{\lambda,n,s}(z)$ there.
Its coefficients are nonnegative, so
\[
 2^s[z^s]F_{\lambda,n,s}(z)\le F_{\lambda,n,s}(2)
 =\exp\left(4b_2(n,s)+\sum_{j\ts\ge\ts 3}2^jb_j(n,s)\right)
 \le M(B)2^{s/2}.
\]
Here $M(B)>0$ is a constant depending only on $B$. The preceding
estimates give the last inequality for all large $n$.
Dividing by $2^s$ gives the
uniform summable bound
\begin{equation}\label{eq:geometric-support}
 T_{\lambda,s}\le M(B)2^{-s/2}
          \qquad(1\le s\le\varepsilon\sqrt n).
\end{equation}

For each fixed $s$, every type containing a transposition or a cycle
of length at least four contributes $o(1)$ uniformly: the corresponding
$b_2$ or $b_j$ tends to zero. The only remaining types are $3^q$.
Lemmas~\ref{lem:contents} and~\ref{lem:three-cycles} give, uniformly,
\[
 r_\lambda(3^q)
    =\frac{\Sigma_\lambda(3)^q}{(n)_{3q}}+O_{B,q}(n^{-q-1})
    =o(n^{-q}).
\]
Their classes have $O_q(n^{3q})$ elements, so their cubic contributions
also tend to zero. Thus, for each fixed $s$,
\[
 \max_{\lambda\ts\in\ts \HH_A(n)}T_{\lambda,s}\longrightarrow0.
\]
Combining this with \eqref{eq:geometric-support},
first truncate at a fixed support and then let the truncation grow, to
obtain
\begin{equation}\label{eq:small-support}
 \max_{\lambda\ts\in\ts \HH_A(n)}
       \sum_{1\ts\le\ts s(\pi)\ts\le\ts \varepsilon\sqrt n}|r_\lambda(\pi)|^3
       \longrightarrow0.
\end{equation}

For larger support, \eqref{eq:orthogonality} and
Lemma~\ref{lem:class-size} imply
\begin{align*}
 \sum_{s(\pi)\ts>\ts \varepsilon\sqrt n}|r_\lambda(\pi)|^3
 &\le e^{3\Delta(\lambda)}
          \max_{s(\pi)\ts>\ts \varepsilon\sqrt n}|C_\pi|^{-1/2}\\
 &\le e^{3A\sqrt n}
          (e/\sqrt n)^{\varepsilon\sqrt n/2}\longrightarrow0.
\end{align*}
This estimate is uniform, completing the proof.
\end{proof}

\subsection{Proofs of the Kronecker results}
\begin{proof}[Proof of Theorem~\ref{thm:uniform}]
By \eqref{eq:character-sum} and H\"older's inequality,
\[
 |R_n(\lambda,\mu,\nu)-1|
 \le\left(\sum_{\pi\ts\ne\ts e}|r_\lambda(\pi)|^3\right)^{1/3}
     \left(\sum_{\pi\ts\ne\ts e}|r_\mu(\pi)|^3\right)^{1/3}
     \left(\sum_{\pi\ts\ne\ts e}|r_\nu(\pi)|^3\right)^{1/3}.
\]
Proposition~\ref{prop:cubic} makes the right side $o(1)$ uniformly.
It is eventually less than $1$, which also proves positivity.
\end{proof}

\begin{proof}[Proof of Corollary~\ref{cor:ppy310}]
Choose a common $A$ for the three sequences and apply
Theorem~\ref{thm:uniform}. For all sufficiently large $n$,
\[
 g(\lambda^{(n)},\mu^{(n)},\nu^{(n)})
       \ge\tfrac12\sqrt{n!}\,e^{-3A\sqrt n}.
\]
For the upper bound, abbreviate
$g:=g(\lambda^{(n)},\mu^{(n)},\nu^{(n)})$ and omit the superscripts
on the dimensions. Symmetry of $g$ and comparison of dimensions in
tensor products give
\[
 g f^\nu\le f^\lambda f^\mu,
 \qquad g f^\mu\le f^\lambda f^\nu.
\]
Multiplying these inequalities and canceling $f^\mu f^\nu>0$ yields
$g^2\le(f^\lambda)^2$. Thus $g\le f^\lambda\le\sqrt{n!}$,
where the last inequality follows from
$\sum_{\alpha\ts\vdash\ts n}(f^\alpha)^2=n!$.
Together with the lower bound and the relative asymptotic, this proves
\eqref{eq:ppy310}. The same upper bound holds for every triple, so
$K(n)\le\sqrt{n!}$; dividing the displayed lower bound by $K(n)$
also gives the comparison with the largest coefficient stated after
the corollary.
\end{proof}

\section{Staircase and square Kronecker coefficients}\label{sec:shapes}

We prove Theorem~\ref{thm:shapes} by summing the character formula over
cycle types. We extract the contribution of disjoint $3$-cycles and
show that all other types have total contribution tending to zero.
For staircases, the uniform estimate uses Murnaghan--Nakayama and a
dimension-ratio bound; for squares, an exterior-power model gives
extra decay for even cycles. We finish by computing the dimension
scale from the hook-length formula.

Recall $\rho_k:=(k-1,k-2,\ldots,1)$ and $\delta_k:=(k^k)$.
Throughout this section, $n$ denotes the size of the diagram under
discussion, so it runs through triangular numbers for $\rho_k$ and
squares for $\delta_k$. Both families are $2$-balanced, meaning that
their first row and first column have length at most $2\sqrt n$.

\subsection{Cycle types and the contribution of three-cycles}
Write a cycle type without fixed points as
$\eta:=(2^{m_2}3^{m_3}\cdots)$, and put
\[
 s:=\sum_{j\ts\ge\ts 2}jm_j,\qquad c:=\sum_{j\ts\ge\ts 2}m_j,\qquad
 b:=\sum_{\substack{j\ts\ge\ts 2\\j\text{ even}}}m_j,\qquad
 z_\eta:=\prod_{j\ts\ge\ts 2}j^{m_j}m_j!.
\]
For $\lambda\vdash n$, abbreviate
$r_\lambda(\eta):=r_\lambda(\eta1^{n-s})$, and define
\begin{equation}\label{sh:class-contribution}
 A_\lambda(\eta):=\frac{(n)_s}{z_\eta}r_\lambda(\eta)^3
       \quad(s\le n),\qquad A_\lambda(\eta):=0\quad(s>n).
\end{equation}
The empty type has contribution one. Counting conjugacy classes in
\eqref{eq:character-sum} gives
\begin{equation}\label{sh:class-sum}
 R_n(\lambda,\lambda,\lambda)=\sum_\eta A_\lambda(\eta),
\end{equation}
where the sum includes all partitions with parts at least two and the
empty partition. The terms can have either sign.

We retain the normalized-character conventions
$r_\lambda(\eta):=\chi^\lambda(\eta1^{n-s})/f^\lambda$ and
$\Sigma_\lambda(\eta):=(n)_s r_\lambda(\eta)$ from
Section~\ref{subsec:characters}, where $(n)_s:=n(n-1)\cdots(n-s+1)$.

\begin{lemma}\label{lem:shape-threes}
For staircases and squares of size $n\ge3$, one has
\[
 r_{\rho_k}(3)=-\frac1{2(n-2)},\qquad
 r_{\delta_k}(3)=-\frac1{n-2},
\]
with $n$ interpreted separately in the two formulas. Put $a:=1/2$ for
the staircase and $a:=1$ for the square. For every fixed $q\ge0$,
\begin{equation}\label{sh:fixed-threes}
 r_\lambda(3^q)=(-a)^q n^{-q}+O_q(n^{-q-1}),
 \qquad
 A_\lambda(3^q)\longrightarrow\frac{(-a^3/3)^q}{q!}.
\end{equation}
\end{lemma}
\begin{proof}
Recall that the $j$th content moment is
\[
 M_j(\lambda):=\sum_{b\ts\in\ts \lambda}c(b)^j
       =\sum_{(i,h)\ts\in\ts \lambda}(h-i)^j.
\]
For the staircase, put $r:=k-1$ and $n:=r(r+1)/2$. Summing contents
along the diagonals $i+j=t$ gives
\[
 M_2(\rho_k)
  =\sum_{t\ts=\ts 2}^{r+1}\sum_{i\ts=\ts 1}^{t-1}(t-2i)^2
  =\frac13\sum_{t\ts=\ts 2}^{r+1}t(t-1)(t-2)
  =\frac{n(n-1)}3.
\]
For the square, $n=k^2$ and
\[
 M_2(\delta_k)
  =2k\sum_{i\ts=\ts 1}^k i^2-2\left(\sum_{i\ts=\ts 1}^k i\right)^2
  =\frac{n(n-1)}6.
\]
Equation~\eqref{eq:content-character} therefore gives
$\Sigma_\lambda(3)=-a n(n-1)$ and the single-cycle formulas.
The square formula is also a special case of
\cite[Corollary~1.3]{MS22}.
Both families are $B$-balanced for a fixed $B$. Apply
Lemma~\ref{lem:three-cycles} and divide by $(n)_{3q}$ to obtain the
first assertion of \eqref{sh:fixed-threes}. The class size
$(n)_{3q}/(3^q q!)$ then gives the second.
\end{proof}

We will prove, separately for each family, a uniform bound of the form
\begin{equation}\label{sh:uniform-template}
 |r_\lambda(\eta)|\le C^s n^{-(s-c)/2}\qquad(s\le n),
\end{equation}
with $C$ independent of $n$ and $\eta$. It implies
\begin{equation}\label{sh:cycle-majorant}
 |A_\lambda(\eta)|\le\prod_{j\ts\ge\ts 2}\frac{y_j^{m_j}}{m_j!},
 \qquad y_j:=\frac{C^{3j}}j n^{(3-j)/2}.
\end{equation}
In particular, $|A_\lambda(3^q)|\le(C^9/3)^q/q!$.
Together with Lemma~\ref{lem:shape-threes}, this summable bound gives
\begin{equation}\label{sh:three-series}
 \sum_{q\ts\ge\ts 0}A_\lambda(3^q)\longrightarrow e^{-a^3/3}.
\end{equation}
Indeed, for a fixed cutoff the convergence follows term by term, and
the remaining absolute tail is bounded by the tail of the exponential
series, independently of $n$.

For later use, if $x_j\ge0$ and $\sum_{j\ts\ge\ts 2}x_j<\infty$, expansion
of the exponential series gives
\begin{equation}\label{sh:exponential-sum}
 \sum_\eta\prod_{j\ts\ge\ts 2}\frac{x_j^{m_j}}{m_j!}
       =\exp\left(\sum_{j\ts\ge\ts 2}x_j\right).
\end{equation}
The sum restricted to types other than $3^q$ is
$e^{x_3}(\exp(\sum_{j\ts\ne\ts 3}x_j)-1)$. These are identities for
nonnegative upper bounds. In particular, for $t>1$,
\begin{equation}\label{sh:weighted-tail}
 \sum_{|\eta|\ts>\ts L}\prod_{j\ts\ge\ts 2}\frac{x_j^{m_j}}{m_j!}
       \le t^{-L}\exp\left(\sum_{j\ts\ge\ts 2}t^jx_j\right)
\end{equation}
whenever the series on the right converges. This follows by
multiplying each summand by the upper bound $t^{|\eta|-L}$.

\subsection{A uniform character bound for staircases}
A \defn{border strip}, or \defn{rim hook}, is a nonempty edge-connected
skew diagram containing no $2$ by $2$ square. Its height is one less
than its number of occupied rows. For an integer $p\ge2$, a partition
is a \defn{$p$-core} if it has no removable rim hook of length $p$.
For the beta-set description of cores and rim-hook removal, see
\cite[Chapter~I, \S1, Example~8]{Mac95}.
The Murnaghan--Nakayama rule states
that, for an $L$-cycle $\tau$ disjoint from $\sigma$,
\begin{equation}\label{sh:MN}
 \chi^\mu(\sigma\tau)
  =\sum_{\substack{\nu\ts\subseteq\ts \mu\\
                  \mu/\nu\text{ a border strip of size }L}}
          (-1)^{\operatorname{ht}(\mu/\nu)}\chi^\nu(\sigma);
\end{equation}
see \cite[Theorem~4.10.2]{Sag01}. The permutation $\sigma$ is regarded as acting on
the remaining letters.

For a partition $\mu$ with at most $r$ rows, pad it with zeros and
define its beta-set by
\[
 \mathcal B_r(\mu):=\{\mu_i+r-i:1\le i\le r\}.
\]
Removing a border strip of length $L$ replaces one $b$ in this set by
$b-L\ge0$, which must be absent from the original set. To check this,
suppose the strip occupies rows $p,\ldots,q$. Consecutive removed
intervals overlap in exactly one column, giving
$\nu_i=\mu_{i+1}-1$ for $p\le i<q$. If $b_i:=\mu_i+r-i$, the new
beta-set replaces $b_p$ by $\nu_q+r-q=b_p-L$. Each choice of $b_p$
determines the new beta-set and hence $\nu$, so there are at most $r$
possible removals of any prescribed length.

For $\rho_k=(r,r-1,\ldots,1)$, the beta-set is
$\{1,3,\ldots,2r-1\}$. An even $L>0$ makes $b-L$ either negative
or already present. Hence \eqref{sh:MN}, with an even cycle removed
first, proves
\begin{equation}\label{sh:stair-even-zero}
 r_{\rho_k}(\eta)=0\quad\text{if $\eta$ has an even part}.
\end{equation}

In particular, no rim hook of length two can be removed from
$\rho_k$, so the staircase is a $2$-core. The beta-set argument
excludes rim hooks of every positive even length, which is the
property used in \eqref{sh:stair-even-zero}.

\begin{lemma}\label{lem:stair-dimension-ratio}
Let $\lambda:=\rho_k$, $r:=k-1$, and $n:=r(r+1)/2$. For every
$\nu\subseteq\lambda$ of size $n-s$,
\begin{equation}\label{sh:stair-dimension-ratio}
 \frac{f^\nu}{f^\lambda}
   \le\frac{[e^{1/2}(r+1)]^s}{(n)_s}
   \le\left(\frac{e^{3/2}(r+1)}n\right)^s.
\end{equation}
Consequently, \eqref{sh:uniform-template} holds for staircases with
$C=16$.
\end{lemma}
\begin{proof}
Pad $\nu$ to length $r$, and put
\[
 b_i:=2(r-i)+1,\qquad a_i:=\nu_i+r-i,\qquad
 u_i:=b_i-a_i=\lambda_i-\nu_i.
\]
Then $u_i\ge0$ and $\sum_i u_i=s$. The determinantal dimension
formula \cite[Theorem~3.11.1]{Sag01}, in its Vandermonde form
\cite[Chapter~I, \S7, Example~6]{Mac95}, gives
\[
 \frac{f^\nu}{f^\lambda}
  =\frac1{(n)_s}\prod_i\frac{b_i!}{a_i!}
                \prod_{i\ts<\ts j}\frac{a_i-a_j}{b_i-b_j}.
\]
We have $b_i!/a_i!\le b_i^{u_i}$ and
$(a_i-a_j)/(b_i-b_j)=1+(u_j-u_i)/(2(j-i))>0$.
Writing $H_m:=\sum_{t\ts=\ts 1}^m1/t$, with $H_0:=0$, and using
$\log(1+x)\le x$, we obtain
\[
 \log\prod_{i\ts<\ts j}\frac{a_i-a_j}{b_i-b_j}
   \le\frac12\sum_i u_i(H_{i-1}-H_{r-i}).
\]
The elementary bounds $H_{i-1}\le1+\log i$ and
$H_{r-i}\ge\log(r-i+1)$ imply
\begin{align*}
 b_i e^{(H_{i-1}-H_{r-i})/2}
 &\le e^{1/2}(2r-2i+1)\sqrt{\frac{i}{r-i+1}}\\
 &\le2e^{1/2}\sqrt{i(r-i+1)}\le e^{1/2}(r+1).
\end{align*}
This proves the first inequality in
\eqref{sh:stair-dimension-ratio}; $(n)_s\ge(n/e)^s$ proves the second.

To bound a character, remove its $c$ nontrivial cycles first in
\eqref{sh:MN}. Every intermediate diagram has at most $r$ rows, so
there are at most $r^c$ removal chains. Each terminal character is the
dimension of a subdiagram with $n-s$ cells. Thus
\[
 |r_\lambda(\eta)|
  \le r^c\left(\frac{e^{3/2}(r+1)}n\right)^s
  \le16^s n^{-(s-c)/2},
\]
using $r,r+1\le2\sqrt n$ and $c\le s/2$.
\end{proof}

\begin{proof}[Proof of \eqref{eq:staircase-asymptotic}]
Lemma~\ref{lem:stair-dimension-ratio} justifies
\eqref{sh:three-series} with $a=1/2$. By
\eqref{sh:stair-even-zero}, all remaining contributing types have
odd parts at least three and at least one part at least five. With
$C=16$, the parameters in \eqref{sh:cycle-majorant} satisfy
\[
 y_3=C^9/3,\qquad
 U_n:=\sum_{\substack{j\ts\ge\ts 5\\j\text{ odd}}}y_j
   \le\frac{C^{15}}{5n(1-C^6/n)}=O(n^{-1}).
\]
Equation~\eqref{sh:exponential-sum} bounds the total absolute
contribution of these remaining types by
$e^{y_3}(e^{U_n}-1)=O(n^{-1})$. Equations~\eqref{sh:class-sum} and
\eqref{sh:three-series} now give the staircase formula.
\end{proof}

\subsection{An exterior-power formula for square characters}
Let $k\ge2$, and put $\lambda:=\delta_k$, $n:=k^2$, and $M:=2k-1$.
On $U:=\CC^{2k}$ with orthonormal basis $u_0,\ldots,u_M$, define
\[
 Fu_i:=\sqrt{i(M-i+1)}\,u_{i-1},\qquad
 E:=F^*,\qquad Hu_i:=(2i-M)u_i,
\]
with out-of-range vectors interpreted as zero. The identities
$[E,F]=H$ and $[E,H]=-2E$ give
\begin{equation}\label{sh:matrix-A}
 A:=e^EFe^{-E}=F+H-E,
 \qquad \|A\|\le4k,\qquad A^{2k}=0.
\end{equation}
Indeed, $\|F\|=\|E\|\le k$ and $\|H\|<2k$. In particular, $A$ is
tridiagonal.

Give $\mathcal W:=\bigwedge^kU$ its induced Hermitian inner product.
For an operator $B$ on $U$, its additive action is
\[
 \dGamma(B)(w_1\wedge\cdots\wedge w_k)
   :=\sum_{i\ts=\ts 1}^k w_1\wedge\cdots\wedge Bw_i\wedge\cdots\wedge w_k.
\]
The tensor-product description immediately gives
\begin{equation}\label{sh:exterior-identities}
 \|\dGamma(B)\|\le k\|B\|,\qquad
 e^{\dGamma(B)}=\bigwedge\nolimits^k(e^B),\qquad
 [\dGamma(B),\dGamma(C)]=\dGamma([B,C]).
\end{equation}
Conjugation by $\bigwedge^kG$ sends $\dGamma(B)$ to
$\dGamma(GBG^{-1})$. Put
\[
 v_+:=u_k\wedge\cdots\wedge u_{2k-1},\qquad
 v_-:=u_0\wedge\cdots\wedge u_{k-1},\qquad
 B_j:=\dGamma(A^j).
\]
The operators $B_j$ commute.

\begin{lemma}\label{lem:square-matrix-formula}
For a type $\eta=(2^{m_2}3^{m_3}\cdots)$ of size $s\le n$,
\begin{equation}\label{sh:square-moment}
 (n)_s r_\lambda(\eta)
    =\left\langle v_+,\left(\prod_{j\ts\ge\ts 2}B_j^{m_j}\right)v_+\right\rangle.
\end{equation}
\end{lemma}
\begin{proof}
In power-sum variables $p_1,p_2,\ldots$, define the complete symmetric
functions by
\[
 h_m:=[t^m]\exp\left(\sum_{j\ts\ge\ts 1}\frac{p_jt^j}{j}\right)
       \quad(m\ge0),
 \qquad h_m:=0\quad(m<0).
\]
Set $\gamma_i:=\sqrt{i!/(M-i)!}$. Since
$F^{b-a}u_b=(\gamma_b/\gamma_a)u_a$ for $a\le b$, the matrix entries satisfy
\[
 \left[\exp\left(\sum_{j\ts\ge\ts 1}\frac{p_jF^j}{j}\right)\right]_{a,b}
    =\frac{\gamma_b}{\gamma_a}h_{b-a}\qquad(0\le a,b\le M).
\]
The matrix element of its exterior power from $v_+$ to $v_-$ is the
minor on rows $0,\ldots,k-1$ and columns $k,\ldots,2k-1$.
Factoring the $\gamma$'s and applying Jacobi--Trudi gives
\begin{equation}\label{sh:square-Schur}
 \left\langle v_-,
    \exp\left(\sum_{j\ts\ge\ts 1}\frac{p_j\dGamma(F^j)}j\right)v_+\right\rangle
       =\mathcal H_k s_{(k^k)}(p),
 \qquad
 \mathcal H_k:=\prod_{i\ts=\ts 0}^{k-1}\frac{(k+i)!}{i!}
             =\frac{n!}{f^\lambda}.
\end{equation}
The remaining determinant is $\det(h_{k+j-i})_{i,j\ts=\ts 1}^k$.
The scalar factored out is
\[
 \frac{\prod_{b\ts=\ts k}^{2k-1}\gamma_b}
      {\prod_{a\ts=\ts 0}^{k-1}\gamma_a}=\mathcal H_k.
\]
The last equality is the hook-length formula for the square.

Apply $\prod_{j\ts\ge\ts 2}(j\partial_{p_j})^{m_j}$ to
\eqref{sh:square-Schur}, and then evaluate at $p_1=1$ and $p_j=0$ for
$j\ge2$. The Frobenius character formula
\[
 s_\lambda(p)=\sum_{\omega\ts\vdash\ts n}
                    \frac{\chi^\lambda(\omega)}{z_\omega}p_\omega
\]
shows that only the type $\eta1^{n-s}$ survives. Differentiation
cancels its factors $j^{m_j}m_j!$, so the right side becomes
$\mathcal H_k\chi^\lambda(\eta1^{n-s})/(n-s)!=(n)_sr_\lambda(\eta)$.
Here Jacobi--Trudi is \cite[Theorem~4.5.1]{Sag01}, and the
Frobenius character formula is \cite[Theorem~4.6.4, (4.23)]{Sag01}. Writing $\widehat F:=\dGamma(F)$, the left side gives
\begin{equation}\label{sh:moment-before-conjugation}
 (n)_s r_\lambda(\eta)
   =\left\langle v_-,e^{\widehat F}
               \left(\prod_{j\ts\ge\ts 2}\dGamma(F^j)^{m_j}\right)v_+\right\rangle.
\end{equation}

To change the left vector, put $w:=e^Fe^{-E}e^F$. Directly,
\[
 wu_i=(-1)^{M-i}u_{M-i}.
\]
One verification identifies $u_i$ with
$\sqrt{\binom Mi}x^{M-i}y^i$: then $F=x\partial_y$,
$E=y\partial_x$, and $w$ sends $x$ to $-y$ and $y$ to $x$.
Thus $w$ is unitary. For $W:=\bigwedge^kw$, the signs from this formula
and from reversing the wedge cancel, giving $Wv_+=v_-$.
Also $\widehat E v_+=0$ for $\widehat E:=\dGamma(E)$, since every
raised occupied index either repeats another one or leaves the range.
Hence $v_+^*\widehat F=0$, and
\[
 v_-^*e^{\widehat F}
   =v_+^*W^{-1}e^{\widehat F}
   =v_+^*e^{-\widehat F}e^{\widehat E}
   =v_+^*e^{\widehat E},
 \qquad e^{-\widehat E}v_+=v_+.
\]
Conjugating the product in \eqref{sh:moment-before-conjugation} by
$e^{\widehat E}$ now replaces each $\dGamma(F^j)$ by
$\dGamma(e^EF^je^{-E})=\dGamma(A^j)=B_j$, proving the lemma.
\end{proof}

\subsection{A decay factor for each even cycle}
For a $k$-element subset $I\subseteq\{0,\ldots,2k-1\}$, let $u_I$
denote its increasing wedge. Put $T:=\{k,\ldots,2k-1\}$ and define
\[
 \mathcal E(I):=\sum_{i\ts\in\ts T}i-\sum_{i\ts\in\ts I}i.
\]
This quantity is nonnegative. Let $P_S$ project onto the span of the
$u_I$ with $\mathcal E(I)\le S$. If
$q:=|I\setminus T|=|T\setminus I|$, then
\begin{equation}\label{sh:energy-defect}
 \mathcal E(I)\ge q^2.
\end{equation}
Indeed, the smallest loss when replacing $q$ upper indices by $q$
lower ones is
$\sum_{i\ts=\ts k}^{k+q-1}i-\sum_{i\ts=\ts k-q}^{k-1}i=q^2$.

\begin{lemma}\label{lem:square-centered-operators}
Put $b_j:=\langle v_+,B_jv_+\rangle$. For $j\ge2$ and $S\ge1$,
\begin{equation}\label{sh:centered-operator-bound}
 \|P_S(B_j-b_j\mathrm{Id})P_S\|\le(16k)^j\sqrt S,
 \qquad \|B_j\|\le k(4k)^j.
\end{equation}
Moreover, $b_j=0$ whenever $j$ is even.
\end{lemma}
\begin{proof}
First let $B$ be any matrix of bandwidth $j$ whose entries have
absolute value at most $L$. Put $b_T:=\sum_{i\ts\in\ts T}B_{ii}$.
At $u_I$, the diagonal entry of $\dGamma(B)-b_T\mathrm{Id}$ has
absolute value at most $2qL$. Each off-diagonal entry replaces an
occupied $i\in I$ by an unoccupied $\ell\notin I$ with
$|i-\ell|\le j$. There are at most
\[
 4jq+\frac{j(j+1)}2
\]
such replacements. At most $2jq$ start at one of the $q$ occupied
lower indices, and at most $2jq$ end at one of the $q$ missing upper
indices. All other replacements cross from an upper index to a lower
one, giving at most $1+\cdots+j$ pairs.

By \eqref{sh:energy-defect}, each absolute column sum of
$P_S(\dGamma(B)-b_T\mathrm{Id})P_S$ is therefore at most
\[
 \left((4j+2)\sqrt S+\frac{j(j+1)}2\right)L.
\]
Applying the same argument to the transpose bounds the row sums by
the same number, which also bounds the operator norm. Now take
$B=A^j$. It has bandwidth $j$, and its entries are bounded by
$(4k)^j$. For $j\ge2$ and $S\ge1$, the parenthesis is at most
$4^j\sqrt S$, proving the centered bound. The uncentered bound
follows from \eqref{sh:matrix-A} and \eqref{sh:exterior-identities}.

Finally, the transpose--sign identity
$\chi^{\lambda'}(\pi)=\sgn(\pi)\chi^\lambda(\pi)$
\cite[Chapter~I, \S7, Example~2]{Mac95} makes the character
of a self-conjugate diagram vanish
on odd permutations. For a single even cycle, Lemma~\ref{lem:square-matrix-formula}
therefore gives $b_j=0$. This covers $j\le n$; when $j>n$, the same
assertion follows from $A^{2k}=0$ and $2k\le n$.
\end{proof}

Recall from \eqref{sh:class-contribution} that $\eta=(2^{m_2}3^{m_3}\cdots)$
has support $s:=\sum_{j\ts\ge\ts 2}jm_j$, $c:=\sum_{j\ts\ge\ts 2}m_j$ cycles,
and $b:=\sum_{j\text{ even}}m_j$ even cycles;
$r_\lambda(\eta):=\chi^\lambda(\eta1^{n-s})/f^\lambda$.

\begin{prop}\label{prop:square-uniform-character}
For the square and every nonempty type of support $s\le n$,
\begin{equation}\label{sh:square-uniform-character}
 |r_\lambda(\eta)|\le64^s n^{-(s-c)/2}(s/n)^{b/2}.
\end{equation}
In particular, \eqref{sh:uniform-template} holds with $C=64$.
\end{prop}
\begin{proof}
Since $A^j$ has bandwidth $j$, the operator $B_j$ changes
$\mathcal E(I)$ by at most $j$. Starting at $v_+$, any intermediate
basis vector in the product in \eqref{sh:square-moment} thus has
energy at most $s$, the sum of all cycle lengths. We may insert $P_s$
between its factors without changing the matrix element.
Use the uncentered bound of \eqref{sh:centered-operator-bound} for
odd $j$ and its centered bound, with $b_j=0$, for even $j$. This gives
\begin{align*}
 |(n)_s r_\lambda(\eta)|
 &\le\prod_{\substack{j\ts\ge\ts 2\\j\text{ odd}}}
                [k(4k)^j]^{m_j}
       \prod_{\substack{j\ts\ge\ts 2\\j\text{ even}}}
                [(16k)^j\sqrt s]^{m_j}\\
 &\le16^s k^{s+c-b}s^{b/2}.
\end{align*}
Divide by $(n)_s\ge(n/e)^s$ and use $k=\sqrt n$ and $16e<64$.
\end{proof}

\subsection{Summing the square character contributions}
\begin{proof}[Proof of \eqref{eq:square-asymptotic}]
Proposition~\ref{prop:square-uniform-character} gives
\eqref{sh:three-series} with $a=1$. It remains to show that the total
absolute contribution of all other types tends to zero. Put
$L:=\lfloor n^{3/5}\rfloor$ and $C:=64$.

For supports exceeding $L$, drop the even-cycle factor in
\eqref{sh:square-uniform-character}. Equations~\eqref{sh:cycle-majorant}
and \eqref{sh:weighted-tail}, with $t=2$, give
\[
 \sum_{|\eta|\ts>\ts L}|A_\lambda(\eta)|
   \le2^{-L}\exp\left(\sum_{j\ts\ge\ts 2}
                 \frac{(2C^3)^j}{j}n^{(3-j)/2}\right)
   =\exp\{-L\log2+O(\sqrt n)\}\longrightarrow0.
\]
Here the $j=2$ term in the exponent is $O(\sqrt n)$, the $j=3$
term is constant, and the remaining terms sum to $O(n^{-1/2})$.

For $s\le L$, retain the even-cycle factor. Cubing
\eqref{sh:square-uniform-character} and multiplying by the class
size gives
\[
 |A_\lambda(\eta)|
   \le\prod_{j\ts\ge\ts 2}\frac{x_j^{m_j}}{m_j!},
 \qquad
 x_j:=\begin{cases}
        C^{3j}n^{(3-j)/2}/j,&j\text{ odd},\\[2pt]
        C^{3j}L^{3/2}n^{-j/2}/j,&j\text{ even}.
      \end{cases}
\]
Then $x_3=C^9/3$ is constant, while geometric-series estimates give
\[
 \sum_{\substack{j\ts\ge\ts 5\\j\text{ odd}}}x_j=O(n^{-1}),
 \qquad
 \sum_{\substack{j\ts\ge\ts 2\\j\text{ even}}}x_j
       =O(L^{3/2}/n)=O(n^{-1/10}).
\]
In \eqref{sh:exponential-sum}, sum over all types other than $3^q$;
dropping the support restriction only increases this nonnegative
upper bound. Their small-support contribution is at most
\[
 e^{x_3}\left(\exp\left(\sum_{j\ts\ne\ts 3}x_j\right)-1\right)
       =O(n^{-1/10}).
\]
Together with the large-support estimate and
\eqref{sh:three-series}, this proves the square formula in
Theorem~\ref{thm:shapes}.
\end{proof}

\subsection{The dimension scale}
The following calculation explains the dimension range of
Theorem~\ref{thm:shapes} and records the resulting growth rates.
Recall the dimension deficit
$\Delta(\lambda):=\log(\sqrt{n!}/f^\lambda)$ for $\lambda\vdash n$.

\begin{prop}\label{prop:shape-dimensions}
With $n$ denoting the size of each diagram separately,
\begin{align}
 \Delta(\rho_k)
    &=\left(\frac32\log2-1\right)n+O(\sqrt n\log n),
       \label{sh:stair-deficit}\\
 \Delta(\delta_k)
    &=(2\log2-1)n+O(\sqrt n\log n).
       \label{sh:square-deficit}
\end{align}
Consequently, for $\lambda_k=\rho_k$ or $\delta_k$,
\[
 g(\lambda_k,\lambda_k,\lambda_k)
   =\sqrt{n!}\exp\{-\gamma n+O(\sqrt n\log n)\},
\]
where $\gamma:=\tfrac92\log2-3$ for staircases and
$\gamma:=6\log2-3$ for squares.
\end{prop}
\begin{proof}
Let $\mathcal H(\lambda):=\prod_{b\ts\in\ts \lambda}h_\lambda(b)$.
The hook-length formula gives
$\Delta(\lambda)=\log\mathcal H(\lambda)-\tfrac12\log n!$.
For the staircase with $r:=k-1$,
\[
 \mathcal H(\rho_k)=\prod_{d\ts=\ts 0}^{r-1}(2d+1)^{r-d}.
\]
Integral comparison, with an $O(r\log r)$ error at the boundary,
gives
\begin{align*}
 \log\mathcal H(\rho_k)
  &=n\log r+r^2\int_0^1(1-x)\log(2x)\,dx+O(r\log r)\\
  &=n\log r+r^2\left(\frac12\log2-\frac34\right)+O(r\log r).
\end{align*}
For the square,
\begin{align*}
 \log\mathcal H(\delta_k)
  &=\sum_{i,j\ts=\ts 0}^{k-1}\log(i+j+1)\\
  &=n\log k+n\int_0^1\int_0^1\log(x+y)\,dx\,dy+O(k\log k)\\
  &=n\log k+n\left(2\log2-\frac32\right)+O(k\log k).
\end{align*}
For these comparisons, monotonicity bounds each sum between adjacent
integral sums; the logarithmic singularity contributes at most the
stated boundary error. Now use $n=r(r+1)/2$ in the first case,
$n=k^2$ in the second, and Stirling's formula. This proves the two
deficit estimates. Finally,
$(f^\lambda)^3/n!=\sqrt{n!}\,e^{-3\Delta(\lambda)}$, and
Theorem~\ref{thm:shapes} supplies the stated growth rates.
\end{proof}

\section{Dimension deficits for LR coefficients}\label{sec:deficits}

We prove Theorem~\ref{thm:deficits} in three steps. A dimension
deficit first bounds the probability that the squared-character
measure preserves a subset of proportional size. We combine this
restriction estimate with Mackey's formula and concentration of the
overlap between the two input blocks. This bounds the input deficits
by the LR coefficient deficit; the restriction-dimension inequality
then controls the output deficit as well.

\subsection{Support and restriction}
For $\rho\vdash r$, recall the normalized restriction norm and
central probability measure from \eqref{eq:Bbeta} and
\eqref{eq:restriction-probability}:
\[
 \beta_\rho(a):=\frac1{\binom ra}
       \sum_{\substack{\alpha\ts\vdash\ts r-a\\\beta\ts\vdash\ts a}}
                  (c^\rho_{\alpha\beta})^2,
 \qquad
 Q_\rho(\pi):=\frac{|\chi^\rho(\pi)|^2}{r!}.
\]
Also recall $\Delta(\rho):=\log(\sqrt{r!}/f^\rho)$ and
$s(\pi):=|\supp\pi|$, the number of moved points.
The measure $Q_\rho$ relates restriction norms to the probability of
preserving a subset.

\begin{lemma}\label{lem:restriction}
Fix $0<\eta\le1/2$, and put
\[
 q_\eta:=\sqrt{\eta^2+(1-\eta)^2},\qquad
 \kappa_\eta:=-\log q_\eta>0.
\]
For $r\ge2$, $\rho\vdash r$ and $\eta r\le a\le(1-\eta)r$,
\begin{equation}\label{eq:restriction-estimate}
 \beta_\rho(a)\le(r+3)
       \exp\{-\kappa_\eta\Delta(\rho)/\log r\}.
\end{equation}
\end{lemma}
\begin{proof}
For $\pi\in S_r$, let $p_a(\pi)$ be the probability that a uniformly
chosen $a$-subset is invariant under $\pi$. If $m_j$ counts its
$j$-cycles, including fixed points, then
\[
 p_a(\pi):=\frac{[z^a]\prod_{j\ts\ge\ts 1}(1+z^j)^{m_j}}{\binom ra}.
\]
Since $Q_\rho$ is central, averaging the indicator of a fixed Young
subgroup over conjugates gives
\begin{equation}\label{eq:central-average}
 \beta_\rho(a)=\EE_{Q_\rho}p_a(\pi).
\end{equation}
Put $\theta:=a/r$ and evaluate the polynomial at
$z=\theta/(1-\theta)$. The binomial probability
$\binom ra\theta^a(1-\theta)^{r-a}$ is a mode and is at least
$1/(r+1)$. Therefore
\[
 p_a(\pi)\le(r+1)\prod_{j\ts\ge\ts 2}
              \bigl(\theta^j+(1-\theta)^j\bigr)^{m_j}
       \le(r+1)q_\eta^{s(\pi)}.
\]
The last inequality uses monotonicity of the $\ell^p$ norms of
$(\theta,1-\theta)$ and $\eta\le\theta\le1-\eta$.

Write $D:=\Delta(\rho)$. The character bound $|\chi^\rho(\pi)|\le f^\rho$
gives $Q_\rho(\pi)\le e^{-2D}$. At most
\[
 \sum_{s\ts=\ts 0}^j\binom rs s!\le\sum_{s\ts=\ts 0}^j r^s\le2r^j
\]
permutations move at most $j$ points. Choose
$j:=\lfloor D/\log r\rfloor\le r/2$. The total $Q_\rho$-mass on
these permutations is at most $2e^{-D}$. On the complement,
$s(\pi)>D/\log r$. Equation~\eqref{eq:central-average} gives
\[
 \beta_\rho(a)\le2e^{-D}+(r+1)e^{-\kappa_\eta D/\log r}
       \le(r+3)e^{-\kappa_\eta D/\log r},
\]
because $\kappa_\eta\le\tfrac12\log2<\log r$.
\end{proof}

\subsection{Concentration of the overlap}
For $k+l=n$, let
\begin{equation}\label{eq:overlap-weights}
 w_a:=\frac{\binom ka\binom la}{\binom nk},
 \qquad a_*:=\frac{kl}{n}.
\end{equation}
Vandermonde's identity \eqref{eq:vandermonde} says that $(w_a)$ is a
probability distribution on the integer indices
$0\le a\le\min(k,l)$. We call precisely these indices
\defn{admissible}; they are the values for which both binomial
coefficients in \eqref{eq:overlap-weights} can be nonzero.
Dividing \eqref{eq:mackey-bound} by $\binom nk$ gives, with
$t:=c^\lambda_{\mu\nu}/\sqrt{\binom nk}$,
\begin{equation}\label{eq:weighted-restriction}
 t^2\le\sum_a w_a\beta_\mu(a)\beta_\nu(a).
\end{equation}

\begin{lemma}\label{lem:overlap}
For $k,l\ge1$ and every admissible integer $0\le a\le\min(k,l)$,
\begin{equation}\label{eq:overlap}
 w_a\le(n+1)e^{-8(a-a_*)^2/n}.
\end{equation}
In particular, for $d>0$,
\begin{equation}\label{eq:overlap-tail}
 \sum_{|a-a_*|\ts>\ts dn}w_a\le(n+1)^2e^{-8d^2n}.
\end{equation}
\end{lemma}
\begin{proof}
Let $h(x):=-x\log x-(1-x)\log(1-x)$, continuously extended to
$[0,1]$. The binomial estimates
\[
 \binom rb\le e^{rh(b/r)},\qquad
 \binom nk\ge\frac{e^{nh(k/n)}}{n+1}
\]
follow from the upper bound $1$ and the modal lower bound $1/(n+1)$
for a binomial probability. Hence $w_a\le(n+1)e^{F(a)}$, where
\[
 F(a):=kh(a/k)+lh(a/l)-nh(k/n).
\]
We have $F(a_*)=F'(a_*)=0$ and, in the interior of the interval,
\[
 F''(a)\le-4/k-4/l\le-16/n.
\]
Taylor's theorem gives \eqref{eq:overlap}; continuity handles the
endpoints. Summing at most $n+1$ terms proves \eqref{eq:overlap-tail}.
\end{proof}

\begin{prop}\label{prop:finite-deficit}
Fix $0<\varepsilon\le1/2$. There is $\gamma_\varepsilon>0$ such that,
for $k,l\ge2$, $k+l=n$, and $k,l\ge\varepsilon n$,
\begin{equation}\label{eq:finite-deficit}
 \left(\frac{c^\lambda_{\mu\nu}}{\sqrt{\binom nk}}\right)^2
 \le2(n+3)^2\exp\left\{-\frac{\gamma_\varepsilon
                     (\Delta(\mu)+\Delta(\nu))}{\log n}\right\}.
\end{equation}
\end{prop}
\begin{proof}
Take $d:=\varepsilon^2/4$ and $\eta:=\varepsilon/2$. If
$|a-a_*|\le dn$, then
\[
 \left|a/k-l/n\right|\le\varepsilon/4,
 \qquad \left|a/l-k/n\right|\le\varepsilon/4.
\]
Thus both $a/k$ and $a/l$ belong to $[\eta,1-\eta]$.
Writing $D:=\Delta(\mu)+\Delta(\nu)$, Lemma~\ref{lem:restriction}
bounds this part of \eqref{eq:weighted-restriction} by
\[
 (k+3)(l+3)
       \exp\left\{-\kappa_\eta
            \left(\frac{\Delta(\mu)}{\log k}
                       +\frac{\Delta(\nu)}{\log l}\right)\right\}
       \le(n+3)^2e^{-\kappa_\eta D/\log n}.
\]
For the other terms use $\beta_\rho(a)\le1$ and
Lemma~\ref{lem:overlap}. Since
$D\le\tfrac12\log(k!l!)\le\tfrac12n\log n$, their contribution is at
most
\[
 (n+1)^2e^{-8d^2n}
       \le(n+1)^2e^{-16d^2D/\log n}.
\]
Adding these bounds proves the result with
$\gamma_\varepsilon:=\min(\kappa_\eta,16d^2)>0$.
\end{proof}

\subsection{Proof of the deficit inequality}
For $c^\lambda_{\mu\nu}>0$, recall from
\eqref{eq:intro-delta} the coefficient deficit
$\delta(\lambda;\mu,\nu):=
\log(\sqrt{\binom nk}/c^\lambda_{\mu\nu})$.

\begin{proof}[Proof of Theorem~\ref{thm:deficits}]
Fix $\varepsilon$. For large $n$, \eqref{eq:VK-upper} gives
\[
 D:=\Delta(\mu)+\Delta(\nu)
       \ge v(\sqrt k+\sqrt l)\ge v\sqrt n.
\]
Consequently,
\[
 \log(2(n+3)^2)\le\frac{\gamma_\varepsilon D}{2\log n}
\]
for all sufficiently large $n$, uniformly over the two partitions.
With $t:=e^{-\delta(\lambda;\mu,\nu)}$, Proposition~\ref{prop:finite-deficit}
now implies
\[
 t^2\le e^{-\gamma_\varepsilon D/(2\log n)},
 \qquad
 D\le\frac4{\gamma_\varepsilon}
                 \delta(\lambda;\mu,\nu)\log n.
\]
By \eqref{eq:dimension-restriction}, the sum of the three deficits is
at most $2D+\delta(\lambda;\mu,\nu)$. This proves
\eqref{eq:deficits} for large $n$ with $a_\varepsilon:=0$ and, for
example, $b_\varepsilon:=8/\gamma_\varepsilon+2$. Only finitely many
sizes and triples remain. Increasing $a_\varepsilon$ covers them all.
\end{proof}

\begin{cor}\label{cor:proportions}
Fix $0<\theta<1$ and let $k_n:=\lfloor\theta n\rfloor$, $l_n:=n-k_n$.
For positive coefficients $c^{\lambda^{(n)}}_{\mu^{(n)}\nu^{(n)}}$ with
these sizes, write
$\delta_n:=\delta(\lambda^{(n)};\mu^{(n)},\nu^{(n)})$. Then
\[
 \Delta(\lambda^{(n)})+\Delta(\mu^{(n)})+\Delta(\nu^{(n)})
       =O_\theta(1+\delta_n\log n).
\]
In particular,
\begin{align*}
 \delta_n=O(n/\log n)&\ \Longrightarrow\
       \Delta(\lambda^{(n)}),\Delta(\mu^{(n)}),\Delta(\nu^{(n)})=O_\theta(n),\\
 \delta_n=o(n/\log n)&\ \Longrightarrow\
       \Delta(\lambda^{(n)}),\Delta(\mu^{(n)}),\Delta(\nu^{(n)})=o(n),\\
 \delta_n=O(\sqrt n)&\ \Longrightarrow\
       \Delta(\lambda^{(n)}),\Delta(\mu^{(n)}),\Delta(\nu^{(n)})
                           =O_\theta(\sqrt n\log n).
\end{align*}
\end{cor}
\begin{proof}
For large $n$, both input sizes are at least $\varepsilon n$, where
$\varepsilon:=\tfrac12\min(\theta,1-\theta)$. Apply
Theorem~\ref{thm:deficits}. Each individual deficit is nonnegative and
bounded by their sum.
\end{proof}

For $\theta=1/2$, exponentiating the bound for $\Delta(\lambda)$
proves \eqref{eq:ppy55}. The same lower bound holds for the normalized
dimensions of each input and their product. If $c=0$, these
exponentiated inequalities hold with right side zero. The first
implication in Corollary~\ref{cor:proportions} is precisely the
fixed-proportion extension discussed in \cite[\S7.7]{PPY19}.

\section{Vanishing with three large dimensions}\label{sec:vanishing}

We prove Theorem~\ref{thm:vanishing} by adding a first row to a
partition of maximal dimension. We choose the row longer than the
sum of the two input widths, forcing the LR coefficient to vanish.
Balancedness keeps its length $O(\sqrt n)$, and an exact hook-length
comparison shows that the new partition still has large dimension.
Here, as before,
$\HH_A(n):=\{\lambda\vdash n:
f^\lambda\ge\sqrt{n!}e^{-A\sqrt n}\}$.

The next lemma computes the dimension cost of adding this row.

\begin{lemma}\label{lem:added-row}
Let $\tau\vdash N$, let $L\ge\tau_1$, and put
$\lambda:=(L,\tau_1,\tau_2,\ldots)\vdash n$, where $n:=N+L$. Then
\begin{equation}\label{eq:added-row}
 \frac{f^\lambda}{\sqrt{n!}}
   =\frac{f^\tau}{\sqrt{N!}}\,
       \frac{\sqrt{n!/N!}}
            {\prod_{j\ts=\ts 1}^{L}(L-j+1+\tau'_j)},
\end{equation}
where $\tau'_j:=0$ for $j>\tau_1$.
\end{lemma}
\begin{proof}
Hooks below the new row are unchanged, and its first-row hooks are
$L-j+1+\tau'_j$. The hook-length formula gives the identity.
\end{proof}

\begin{proof}[Proof of Theorem~\ref{thm:vanishing}]
Recall $a_0:=\pi/\sqrt6$ from \eqref{eq:partition-bound} and
that $B$-balanced means that the first row and column have length
at most $B\sqrt n$ for a partition of size $n$.
Lemma~\ref{lem:balance} gives $B:=B(A)$ with
\begin{equation}\label{eq:input-width}
 \mu_1+\nu_1\le2B\sqrt{n/2}=\sqrt2B\sqrt n.
\end{equation}
For each $N$, choose a maximal-dimensional partition $\tau\vdash N$.
By \eqref{eq:partition-bound}, it belongs to $\HH_{a_0}(N)$, and
Lemma~\ref{lem:balance} gives a universal $B_0$ with
$\tau_1,\ell(\tau)\le B_0\sqrt N$.

Fix $D>\max(\sqrt2B,B_0)$. For large $n$, take
\[
 L:=\lceil D\sqrt n\rceil,
 \qquad N:=n-L\ge n/2,
 \qquad \lambda:=(L,\tau_1,\tau_2,\ldots).
\]
Then $L>\tau_1$, so $\lambda$ is a partition, and
$\lambda_1>\mu_1+\nu_1$. The LR rule implies
\begin{equation}\label{eq:LR-width}
 c^\lambda_{\mu\nu}>0\quad\Longrightarrow\quad
                     \lambda_1\le\mu_1+\nu_1.
\end{equation}
Indeed, if $\mu\subseteq\lambda$, the first row of an LR tableau of
shape $\lambda/\mu$ consists entirely of $1$'s, and there are only
$\nu_1$ such entries in its content. If containment fails, the
coefficient is zero already. Hence our coefficient vanishes.

It remains to check its dimension. Put $C:=D+B_0+1$. Every new hook is
at most $L+\ell(\tau)\le C\sqrt n$. By Lemma~\ref{lem:added-row},
\eqref{eq:partition-bound} and $n!/N!\ge(n/2)^L$,
\begin{align*}
 \frac{f^\lambda}{\sqrt{n!}}
 &\ge e^{-a_0\sqrt N}
                \frac{\sqrt{n!/N!}}{(C\sqrt n)^L}\\
 &\ge e^{-a_0\sqrt n}(\sqrt2C)^{-L}
  \ge\exp\{-[a_0+(D+1)\log(\sqrt2C)]\sqrt n\}.
\end{align*}
The expression in square brackets is a permissible $A'$, depending
only on $A$.
\end{proof}

\begin{remark}\label{rem:zero}
To obtain the sequences in \cite[Conjecture~4.16]{PPY19}, choose
$\mu^{(n)}=\nu^{(n)}$ of maximal dimension in $S_{n/2}$ and apply
Theorem~\ref{thm:vanishing} with $A=a_0$. The input and output
dimensions then all have deficit $O(\sqrt n)$, while the coefficient
is zero for every sufficiently large even $n$. With
$\log0:=-\infty$, the expression in that conjecture is eventually
$+\infty$. If a strictly positive coefficient is required in its
logarithmic formulation, that additional variant is not proved by this
construction. The improvement of \cite[Proposition~4.17]{PPY19}
holds independently of this convention.
\end{remark}

\section{The exact stabilization threshold}\label{sec:stabilization}

We prove Theorem~\ref{thm:stabilization} by characterizing equality
in the dimension bound for an LR coefficient. Row and column
stabilizers give a strict bound unless the relevant skew diagram is
a rook strip. We then determine the smallest outer partition that
can contain such a strip and construct one at every size beyond
the threshold. Recall from \eqref{eq:intro-maxima} that $D(k)$ is
the largest irreducible dimension of $S_k$, $C(n,k)$ is the largest
LR coefficient with input sizes $k,n-k$, and $N_k:=\binom{k+1}{2}$.

\subsection{Equality and rook strips}
Recall that a \defn{rook strip} has at most one cell in each row and
each column. We write $[W:V]$ for the multiplicity of an irreducible
module $V$ in a semisimple module $W$; $\one$ denotes the trivial
character and $\sgn$ the sign character.
For a skew diagram $\eta$ with $k$ cells, its skew Specht module
$S^\eta$ has Frobenius characteristic $s_\eta$
\cite[Chapter~I, \S7, Example~3]{Mac95}; hence
\begin{equation}\label{eq:skew-Specht}
 S^{\lambda/\nu}\cong\bigoplus_{\mu\ts\vdash\ts k}
                    (V^\mu)^{\oplus c^\lambda_{\mu\nu}}.
\end{equation}
We use the row-tabloid and polytabloid constructions, recalled in the
proof below; see \cite[\S\S2.1, 2.3]{Sag01}.

\begin{lemma}\label{lem:row-column}
Let $\eta:=\lambda/\nu$ have $k$ cells, and let $\mu\vdash k$. Then
$c^\lambda_{\mu\nu}\le f^\mu$. For $k\ge2$, let $t\in S_k$ be a
transposition. If a row of $\eta$ has at least two cells, then
\begin{equation}\label{eq:row-bound}
 c^\lambda_{\mu\nu}\le\frac{f^\mu+\chi^\mu(t)}2.
\end{equation}
If a column of $\eta$ has at least two cells, then
\begin{equation}\label{eq:column-bound}
 c^\lambda_{\mu\nu}\le\frac{f^\mu-\chi^\mu(t)}2.
\end{equation}
\end{lemma}
\begin{proof}
Fix a bijective filling $T$ of $\eta$ by $[k]$, with row and column
stabilizers $R_T$ and $K_T$. The row-tabloid module is
$M^\eta:=\Ind_{R_T}^{S_k}\one$. The skew Specht module is its submodule
spanned by the translates of the polytabloid
\[
 e_T:=\sum_{\sigma\ts\in\ts K_T}\sgn(\sigma)\,\sigma\{T\}.
\]
By semisimplicity and Frobenius reciprocity,
\[
 c^\lambda_{\mu\nu}\le[M^\eta:V^\mu]
       =\dim(V^\mu)^{R_T}\le f^\mu.
\]
If a row is repeated, $R_T$ contains a transposition. Its fixed space
on $V^\mu$ has dimension $(f^\mu+\chi^\mu(t))/2$, proving
\eqref{eq:row-bound}.

For the column bound, $K_T$ acts on $e_T$ by its sign character.
The resulting homomorphism
\[
 \Ind_{K_T}^{S_k}\sgn_{K_T}\longrightarrow S^\eta
\]
is surjective, since the translates of $e_T$ span. Thus the
multiplicity is bounded by the dimension of the subspace of $V^\mu$
on which $K_T$ acts by sign. If a column is repeated, that subspace is
contained in the $(-1)$-eigenspace of a transposition, of dimension
$(f^\mu-\chi^\mu(t))/2$.
\end{proof}

\begin{lemma}\label{lem:equality}
Let $\nu\subseteq\lambda$, $k:=|\lambda/\nu|$, and $\mu\vdash k$.
If $f^\mu>1$, then
\begin{equation}\label{eq:rook-equality}
 c^\lambda_{\mu\nu}=f^\mu
       \quad\Longleftrightarrow\quad
       \lambda/\nu\text{ is a rook strip}.
\end{equation}
For a rook strip, equality holds for every $\mu\vdash k$.
\end{lemma}
\begin{proof}
For $f^\mu>1$, one has $|\chi^\mu(t)|<f^\mu$. Otherwise every
eigenvalue of $t$ would be the same sign. Conjugacy then makes all
transpositions act by the same scalar, and they generate $S_k$.
Irreducibility would force dimension one. If $\lambda/\nu$ is not a
rook strip, Lemma~\ref{lem:row-column} consequently gives
\[
 c^\lambda_{\mu\nu}\le\frac{f^\mu+|\chi^\mu(t)|}{2}<f^\mu.
\]
Conversely, for a rook strip both $R_T$ and $K_T$ are trivial. The
tabloid module is regular and each polytabloid equals its tabloid, so
$S^{\lambda/\nu}$ is the regular representation. Each $V^\mu$ occurs
with multiplicity $f^\mu$.
\end{proof}

\subsection{The first occurrence of the maximal value}
\begin{lemma}\label{lem:rook-size}
If $\lambda/\nu$ is a rook strip with $k$ cells, then
$|\lambda|\ge\binom{k+1}{2}$.
\end{lemma}
\begin{proof}
Each occupied row contributes its last cell in $\lambda$. These $k$
cells occupy distinct columns, so the corresponding row lengths of
$\lambda$ are $k$ distinct positive integers. Their sum is at least
$1+\cdots+k$.
\end{proof}

\begin{proof}[Proof of Theorem~\ref{thm:stabilization}]
Lemma~\ref{lem:row-column} gives $C(n,k)\le D(k)$.
Suppose $k\ge3$ and equality holds, and choose a maximizing triple.
Then
\[
 D(k)=c^\lambda_{\mu\nu}\le f^\mu\le D(k).
\]
The shape $(k-1,1)$ has dimension $k-1$, so $D(k)>1$.
Lemma~\ref{lem:equality} forces $\lambda/\nu$ to be a rook strip.
Lemma~\ref{lem:rook-size} now gives $n\ge N_k$.

For sufficiency, let $r:=n-N_k\ge0$ and take
\[
 \lambda:=(k+r,k-1,k-2,\ldots,1),
 \qquad \nu:=(k-1+r,k-2,\ldots,1),
\]
omitting zero parts. Then $|\nu|=n-k$ and $\lambda/\nu$ is a rook
strip with $k$ cells. By Lemma~\ref{lem:equality},
$c^\lambda_{\mu\nu}=f^\mu$ for every $\mu\vdash k$. Choose $\mu$
of maximal dimension. This is the construction in
\cite[Theorem~4.23]{PPY19}.

Finally, $D(0)=D(1)=D(2)=1$. The upper bound still applies, and
$c^{(n)}_{(k),(n-k)}=1$ gives equality for every $n\ge k$, with the
empty partition convention. Thus the first occurrence is $N_k$ for
$k\ne2$ and is $2$ for $k=2$ on the stated domain.
\end{proof}

\end{document}